\documentclass{article}

\usepackage[parfill]{parskip}
\usepackage{amsmath}
\usepackage{amssymb}
\usepackage{mathtools}
\usepackage{amsfonts}
\usepackage{amsthm}
\usepackage{enumitem}
\usepackage{thmtools, thm-restate}
\usepackage{biblatex}
\usepackage{import}
\usepackage{hyperref}
\usepackage{csquotes}
\usepackage{color}

\declaretheorem[numberwithin=section]{theorem}
\declaretheorem[sibling=theorem]{definition}
\declaretheorem[sibling=theorem]{lemma}
\declaretheorem[sibling=theorem]{example}
\declaretheorem[sibling=theorem]{remark}
\declaretheorem[sibling=theorem]{corollary}

\DeclareMathOperator{\Image}{Im}
\DeclareMathOperator*{\argmin}{arg\,min}

\DeclareMathOperator{\sgn}{sgn}

\numberwithin{equation}{section}

\newcommand{\parentheses}[1]{\left(#1\right)}
\newcommand{\abstractUnderlyingSpace}{\mathcal{X}}
\newcommand{\abstractBorelSpace}{\parentheses{\abstractUnderlyingSpace, \tau}}
\newcommand{\probabilityMeasures}[1]{\mathcal{P}\parentheses{#1}}

\newcommand{\density}[2]{\frac{d#1}{d#2}}
\newcommand{\integral}[2]{\int #1 \,d#2}
\newcommand{\LpSpace}[2]{L^{#2}\parentheses{#1}}

\newcommand{\essentiallyBounded}[1]{\mathcal{L}^{\infty}\parentheses{#1}}
\newcommand{\set}[1]{\left\{#1\right\}}

\newcommand{\divergence}[3]{D_{#1}\parentheses{#2 \middle\| #3}}
\newcommand{\reals}{\mathbb{R}}
\newcommand{\extendedReals}{\overline{\reals}}
\newcommand{\nonNegativeSubset}[1]{#1_{\geq 0}}
\newcommand{\integers}{\mathbb{N}}

\newcommand{\othorgonal}[1]{#1^{\perp}}
\newcommand{\closure}[1]{\overline{#1}}

\newcommand{\CDIterate}[2]{#1^{(#2)}}
\newcommand{\CDIntermediateIterate}[3]{#1^{(#2, #3)}}
\newcommand{\absoluteValue}[1]{\left\lvert #1 \right\rvert}
\newcommand{\norm}[1]{\left\lVert #1 \right\rVert}
\newcommand{\specifiedNorm}[2]{\norm{#1}_{#2}}
\newcommand{\sumNorm}[1]{\specifiedNorm{#1}{\ell^2}}
\newcommand{\generalSumNorm}[2]{\specifiedNorm{#1}{\ell_{#2}^2}}
\newcommand{\generalQuotientNorm}[2]{\specifiedNorm{#1}{\ell_{#2}^2 / \sim}}
\newcommand{\lTwoNorm}[1]{\specifiedNorm{#1}{\LpSpace{\mu}{2}}}
\newcommand{\supNorm}[1]{\specifiedNorm{#1}{\infty}}

\newcommand{\openInterval}[2]{\left(#1, #2 \right)}
\newcommand{\closedInterval}[2]{\left[#1, #2\right]}
\newcommand{\innerProduct}[2]{\left\langle #1, #2 \right\rangle}
\newcommand{\lTwoInnerProduct}[2]{\innerProduct{#1}{#2}_{\LpSpace{\mu}{2}}}
\newcommand{\sumInnerProduct}[2]{\innerProduct{#1}{#2}_{\ell^2}}

\newcommand{\gateauxDerivative}[1]{d #1}
\newcommand{\partialGateauxDerivative}[2]{d_{#2} #1}
\newcommand{\factoredDomain}{\hat{H}_\infty}
\newcommand{\quotientNorm}[1]{\specifiedNorm{#1}{\sim}}

\newcommand{\equivalenceClass}[1]{[ #1 ]}

\newcommand{\fixedSizeNorm}[1]{\lVert #1 \rVert}

\newcommand{\friedrich}[2]{\mathcal{C}\parentheses{#1, #2}}
\newcommand{\generalFriedrich}[3]{\mathcal{C}_{#3}\parentheses{#1, #2}}

\newcommand{\X}{\mathcal{X}}

\author{Stephan Eckstein\thanks{Department of Mathematics, University of T\"{u}bingen, stephan.eckstein@uni-tuebingen.de.}, Aziz Lakhal\thanks{Department of Mathematics, University of T\"{u}bingen, aziz.lakhal@uni-tuebingen.de.}}

\title{Exponential convergence of general iterative proportional fitting procedures}

\date{\today}

\begin{document}
\maketitle

\begin{abstract}
	Motivated by the success of Sinkhorn's algorithm for entropic optimal transport, we study convergence properties of iterative proportional fitting procedures (IPFP) used to solve more general information projection problems.
    We establish exponential convergence guarantees for the IPFP whenever the set of probability measures which is projected onto is defined through constraints arising from linear function spaces. This unifies and extends recent results from multi-marginal, adapted and martingale optimal transport. 
    The proofs are based on strong convexity arguments for the dual problem, and the key contribution is to illuminate the role of the geometric interplay between the subspaces defining the constraints. In this regard, we show that the larger the angle (in the sense of Friedrichs) between the linear function spaces, the better the rate of contraction of the IPFP. 
\end{abstract}

\section{Introduction}
Optimization problems over sets of probability measures play a major role in many recent research areas, like variations of optimal transport 
(see, e.g., \cite{backhoff2017causal,beier2022linear,cheridito2025optimal,peyre2019computational,zaev2015monge}), distributionally robust 
optimization (see, e.g., \cite{hu2013kullback, husain2024distributionally, mohajerin2018data, rahimian2019distributionally}), information projections and related problems (see, e.g., \cite{arnese2024convergence,carlier2017convergence,csiszar1975divergence, leonard2012schrodinger}) and many others (see, e.g., \cite{bellec2024optimizing,cheridito2017duality,de2016imprecise,frohlich2024risk,graf2000foundations}).
In many problems, regularization methods are essential for numerical tractability of these optimization problems. The primary example in this regard is entropic regularization for optimal transport (see, e.g., \cite{carlier2017convergence,cuturi2013sinkhorn,nutz2021introduction,pooladian2021entropic}), which allows the problem to be numerically solved via Sinkhorn's algorithm (see, e.g., \cite{cuturi2013sinkhorn, deming1940least, sinkhorn1967diagonal}). The speed and simplicity of this algorithm have been of profound importance, as it made optimal transport easily applicable as a tool for many downstream applications, such as in generative modeling (see, e.g., \cite{de2021diffusion,genevay2018learning}), computational biology (see, e.g., \cite{bunne2022proximal, cang2023screening, demetci2022scot}), image analysis (see, e.g., \cite{ge2021ota,zhou2022rethinking}) and compression (see, e.g., \cite{cui2024sinkd,lei2022neural}), just to name a few. Importantly, entropic regularization turns optimal transport into an information projection problem, and Sinkhorn's algorithm is a particular case of the iterative proportional fitting procedure (IPFP for short). In this paper, we establish exponential convergence guarantees of the IPFP when applied to solving general information projection problems.

\subsection{Motivation and overview of results}
Information projection problems (cf.~\cite{csiszar1975divergence,csiszar2003information}) are of the form
\begin{equation}
\inf_{\pi \in \mathcal{Q}} D_{\rm KL}(\pi \| \theta) \label{eq:IP},
\end{equation}
where $\mathcal{Q} \subset \mathcal{P}(\X)$ is a set of probability measures on a polish space $\X$ and $\theta \in \mathcal{P}(\X)$ is the probability measures which is projected onto this set, while $D_{\rm KL}$ is the Kullback-Leibler divergence, which is the notion of distance used for the projection. While later sections consider more general problems including different divergences, we stick to \eqref{eq:IP} in the introduction for notational simplicity.

The iterative proportional fitting procedure (cf.~\cite{csiszar1975divergence,dykstra1985iterative}), or simply IPFP, is an algorithm to approximately solve the information projection problem \eqref{eq:IP}.\footnote{We should mention that many authors use the term \textsl{IPFP} more narrowly for the algorithm applied to the information projection problem corresponding to entropic optimal transport. We believe, also in line with \cite{csiszar1975divergence}, that the term fits the general procedure perfectly as well.} It based on the idea of writing the set $\mathcal{Q}$ as an intersection $\mathcal{Q} = \cap_{i=1}^N \mathcal{Q}_i$ and iteratively projecting onto the sets in the intersection. That is, one starts with $\pi^{(0)} = \pi^{(0, 0)} = \theta$ and, for $t\geq 0$, sets\footnote{Under mild assumptions which are satisfied later on, the minimizer of \eqref{eq:ipfp} exists uniquely, see also \cite[Theorem 2.1]{csiszar1975divergence}.}
\begin{align}\label{eq:ipfp}
	\begin{split}
\pi^{(t, i)} &:= \argmin_{\pi \in \mathcal{Q}_{i}} D_{\rm KL}(\pi \| \pi^{(t, i-1)}) \text{ for } i=1, \dots, N ~\text{ and } \\
\pi^{(t)} &:= \pi^{(t, 0)} := \pi^{(t-1, N)}.
\end{split}
\end{align}
The idea behind this procedure is that the projection onto $\mathcal{Q}$ may not be directly tractable to compute, but the projections onto the simpler sets $\mathcal{Q}_i$ are tractable and the iterates $\pi^{(t)}$ converge to the projection onto $\mathcal{Q}$.
As an example, in entropic optimal transport the set $\mathcal{Q}$ consists of all joint distributions with fixed marginal distributions $\mu_1$ and $\mu_2$, and Sinkhorn's algorithm in this case uses $N=2$, where $\mathcal{Q}_i$ is the set which requires only one of the marginals $\mu_i$ to be fixed. The projection onto $\mathcal{Q}_i$ admits a closed form solution in this case. Even in situations in which closed form solutions are not available, approximations of the projections within the IPFP can be used (see, e.g., \cite{von2023generalized}).

While qualitative convergence of the IPFP has been studied and is known in many settings for a long time (see, e.g., \cite{csiszar1975divergence,fienberg1970iterative,ruschendorf1995convergence}), the rate of convergence has only more recently been a frequent topic of interest, particularly for the case of entropic optimal transport. In this context, it has been shown that the algorithm converges with a linear rate in fairly general settings (see, e.g., \cite{chen2016entropic,conforti2023quantitative,eckstein2023hilbert}). While this rate of convergence has been established for other optimization problems beyond classical optimal transport as well (see \cite{carlier2022linear,chen2024convergence,eckstein2024computational}), a general treatment for the IPFP for problems of the form \eqref{eq:IP} covering a linear rate of convergence is open.

In this paper, we establish linear rates of convergence of the IPFP for sets of measures which arise from constraints defined by linear function spaces. That is, for linear spaces $H_i \subseteq L^\infty(\X)$, we consider
\[
\mathcal{Q}_i = \left\{\pi \in \mathcal{P}(\X) \,:\, \int h_i \,d\pi = \int h_i \,d\mu \text{ for all }h_i \in H_i\right\}, ~ i=1, \dots, N.
\]
Hereby, $\mu \in \mathcal{P}(\X)$ is a fixed reference measure which is equivalent to $\theta$. We denote by $\closure{H}_i$ the closure of $H_i$ in $L^2(\mu)$. The following specifies the main convergence result of Section \ref{convergenceSection} to the particular setting presented in the introduction.
%

\begin{theorem}\label{thm:intro}
	Assume \eqref{eq:IP} satisfies strong duality and the primal and dual versions of the IPFP coincide (cf.~Section \ref{sec:setting}) and let $\pi^*$ be an optimizer of $\eqref{eq:IP}$. Assume further that
	\begin{enumerate}
		\item the log-densities of $(\pi^{(t)})_{t\in\mathbb{N}}$ and $\pi^*$ with respect to $\mu$ are uniformly bounded by $R > 0$, and 
		\item the sum $\closure{H}_1 + \closure{H}_2 + \dots + \closure{H}_N$ is closed in $L^2(\mu)$.
	\end{enumerate} 
	Then, there exists a constant $C > 0$ and $\rho \in (0, 1)$ such that
	\[
		\big| D_{\rm KL}(\pi^{(t)} \| \mu) - D_{\rm KL}(\pi^* \| \mu) \big| \leq C \rho^t.
	\]
	\begin{proof}
		The statement is a direct consequence of Theorem \ref{exponential_convergence:thm}.
	\end{proof}
\end{theorem}

The proof of the above theorem is based on prior works studying the IPFP via strong convexity methods (mainly \cite{carlier2022linear,eckstein2024computational}). The first assumption, boundedness of the iterates, is fairly standard and we will discuss it later on. 

\paragraph{Relation between the geometry and sum of $H_1, \dots, H_N$.}
The main contribution of this work is to clarify the role of the geometric interplay between the spaces $H_1, H_2, \dots, H_N$, which is expressed through the simple closedness condition of their sum. At first, this condition may appear surprising. However, this seemingly purely topological assumption naturally relates to the topic of angles between spaces, and angles between spaces are more easily seen to be intuitively relevant for projections. Indeed, it is established in \cite{deutsch1995angle} that $\closure{H}_1 + \closure{H}_2$ is closed if and only if the angle (in the sense of Friedrichs \cite{friedrichs1937certain}) between the two spaces is strictly positive (cf.~Definition \ref{friedrichsAngleDefinition} for the precise notion of angle used). More generally, from a dual perspective it is clear that summation of elements from the different spaces is a crucial operation in the IPFP. The closedness of $\closure{H}_1 + \closure{H}_2 + \dots + \closure{H}_N$ precisely yields that the condition number (which is a measure of numerical stability) of the sum operator is finite. To make this more precise, let us write $S(h_1, \dots, h_N) := h_1 + \dots + h_N$ for $h_i \in \closure{H}_i, i=1, \dots, N$, and define the sum operator modulo its kernel by
\begin{align*}
\hat{S} : (\closure{H}_1 \times \dots \times \closure{H}_N)/\ker(S) &\rightarrow \closure{H}_1 + \dots + \closure{H}_N,\\ 
[(h_1, \dots, h_N)] &\mapsto h_1 + \dots + h_N.
\end{align*}
That is, $\hat{S}$ factors out equivalences in $\closure{H}_1 \times \dots \times \closure{H}_N$ which do not alter the sum. While $\closure{H}_1 + \dots + \closure{H}_N \subset L^2(\mu)$ is equipped with  $\|\cdot\|_{L^2(\mu)}$, the space $(\closure{H}_1 \times \dots \times \closure{H}_N)/\ker(S)$ is endowed with the norm 
\[\|[(h_1, \dots, h_N)]\|_\sim^2 := \inf_{r \in \ker(S)} \sum_{i=1}^N \|h_i+r_i\|_{L^2(\mu)}^2.\]

\paragraph{Explicit bounds for the contraction coefficient.} With this notation in place, we can give the form of the contraction coefficient $\rho$ in Theorem \ref{thm:intro}, which is given by
\begin{equation}\label{eq:rhoformula}
\rho = 1 - \frac{1}{N} \left(\frac{\sigma_R}{L_R} \cdot \frac{1}{\|\hat{S}\| \|\hat{S}^{-1}\|}\right)^2.
\end{equation}
The term $\frac{\sigma_R}{L_R}$ solely depends on the boundedness condition in Theorem \ref{thm:intro} and is treated below. Let us instead focus on $\|\hat{S}\| \|\hat{S}^{-1}\|$, which is precisely the condition number of the operator $\hat{S}$. We establish at the start of Section \ref{subsec:results} that this condition number is finite if and only if $\closure{H}_1 + \dots + \closure{H}_N$ is closed, see Remark \ref{rem:iffinvertible}. Further, coming back to the relation between the sum and the angle between the spaces: In Section \ref{general_considerations_subsection} we show that the condition number $\|\hat{S}\| \|\hat{S}^{-1}\|$ can be computed using angles between the subspaces $H_1, \dots, H_N$. 
To be more precise, for $N=2$, Theorem \ref{two_subspaces_bounds:thm} shows that
\[
\left(\|\hat{S}\| \|\hat{S}^{-1}\|\right)^2 = \frac{1 + \friedrich{\closure{H}_1}{ \closure{H}_2}}{1-\friedrich{\closure{H}_1}{ \closure{H}_2}},
\]
where $\friedrich{\closure{H}_1}{ \closure{H}_2} \in [0, 1]$ is the cosine of the angle between the spaces $H_1$ and $H_2$, see Definition \ref{friedrichsAngleDefinition} for details. Again, we emphasize that $\friedrich{\closure{H}_1}{ \closure{H}_2} < 1$ if and only if $\closure{H}_1 + \closure{H}_2$ is closed, which thus makes the condition in Theorem \ref{thm:intro} appear completely natural at last.

For $N>2$, a precise computation of the condition number $\|\hat{S}\| \|\hat{S}^{-1}\|$ through angles is beyond the scope of this paper; however, we can at least bound the condition number through iterative applications of angles between pairs of subspaces, see Theorem \ref{general_subspaces:thm}.

We further emphasize that the formulas involving angles are not just abstractly appealing, but lead to concrete and simple ways to bound contraction coefficients in relevant examples, see Remark \ref{generalizedAngle:remark}. We explore this in the context of the IPFP applied to entropic martingale optimal transport in Section \ref{mot_subsection}, which complements the recent work \cite{chen2024convergence}. At this point, we should also mention an obvious fact: For (multi-marginal) optimal transport, all spaces $H_i$ are orthogonal (see Example \ref{ex:comparisoncarlier}) when $\mu$ is the product measure of the marginals and thus the condition number of $\hat{S}$ is always best possible (equal to one). This may be a reason why the IPFP is such an attractive algorithm in optimal transport settings.

\paragraph{Boundedness condition.} We next discuss the boundedness condition of Theorem \ref{thm:intro}. First, we mention that several recent works, primarily on optimal transport, have treated boundedness (or, more generally, smoothness) conditions of the IPFP (see, e.g., \cite{chiarini2024semiconcavity,chizat2024sharper,eckstein2023hilbert}), and the takeaway message is clear: The better one can bound the iterates of the algorithm, the better the contraction coefficient $\rho$. Our results only modestly contribute to this aspect, as we simply assume boundedness holds and establish the corresponding influence of the bounds on the contraction coefficient $\rho$ in this general framework. This influence is mediated by the fraction $\frac{\sigma_R}{L_R}$ of the the strong convexity constant and Lipschitz constant of the dual penalization corresponding to the divergence used. In case of the relative entropy, this would simply correspond to $\frac{\sigma_R}{L_R} = \exp(2R)$. We believe the pattern of influence of both strong convexity and Lipschitz continuity is primarily interesting in light of more general divergences than $D_{\rm KL}$. As a side note, the effect of the bound $R$ on the contraction coefficient we prove is the same as the one established in \cite{carlier2022linear} for entropic multi-marginal optimal transport, see \mbox{Example \ref{ex:comparisoncarlier}}.

\paragraph{Organization of the paper.} The remainder of the paper is structured as follows: Section \ref{sec:setting} fixes our setting and notation, and introduces the dual perspective of the IPFP together with the standing assumption on strong duality. Notably, this includes generalizations of \eqref{eq:IP} using more general divergences and including a linear cost term. Section \ref{convergenceSection} establishes the general version of Theorem \ref{thm:intro}, that is, the linear convergence of the IPFP. Section \ref{general_considerations_subsection} then explores the relationship between the condition number of the sum operator and angles between subspaces, while Subsection \ref{mot_subsection} showcases these results at the concrete example of entropic martingale transport.

\section{Setting and Notation}\label{sec:setting}

Let $ \abstractBorelSpace $ denote a topological space.
Denote by $ \probabilityMeasures{\abstractUnderlyingSpace} $ the space of probability
measures defined on its respective Borel $ \sigma $-algebra.
Choose a reference measure $ \mu \in \probabilityMeasures{\abstractUnderlyingSpace}$. For 
$ \pi \in \probabilityMeasures{\abstractUnderlyingSpace} $ and a convex function 
$ \varphi : \nonNegativeSubset{\reals} \to \reals $, such that $ \varphi(1) = 0 $, let
$ \divergence{\varphi}{\pi}{\mu} $ denote the $ \varphi $-divergence of $ \pi $ from $ \mu $.

We will now constrain $ \probabilityMeasures{\abstractUnderlyingSpace} $ as follows. Let $ N \geq 1 $ and consider non-trivial (i.e., not equal to $\{0\}$)
\emph{linear} subspaces $ H_i \subset \LpSpace{\mu}{\infty} $, $ i = 1, \ldots, N $. To avoid degenerate cases later on when discussing angles between spaces, we assume that none of the subspaces is completely contained in the sum of the remaining subspaces. With each subspace $ H_i $, we associate
\begin{equation*}
    \mathcal{Q}_i \coloneqq \set{\pi \in \probabilityMeasures{\abstractUnderlyingSpace} : 
    \integral{h_i}{\pi} = \integral{h_i}{\mu} \ \text{ for all } h_i
    \in H_i} \ ,
\end{equation*}
and denote the intersection of the above constraint sets by $ \mathcal{Q} = \cap_{i=1}^N \mathcal{Q}_i $. 

Given $ c \in \essentiallyBounded{\mu} $ and $ \varepsilon > 0 $, our goal is to solve the constrained problem
\begin{equation}
    \label{primalProblem}
    \tag{PP}
    \underset{\pi \in \mathcal{Q}}{\inf} J\parentheses{\pi}, ~ \text{ where }
    J\parentheses{\pi} \coloneqq \integral{c}{\pi} + \varepsilon \divergence{\varphi}{\pi}{\mu}.
\end{equation}
We can certainly assume that $ \varepsilon = 1 $, for if not, we replace $ c $ by $ c / \varepsilon $. We hence work with $\varepsilon = 1$ in the remainder of the paper. 

To clarify the relation of \eqref{primalProblem} to the standard information projection \eqref{eq:IP} from the introduction with objective $\tilde{J}(\pi) = D_{\rm KL}(\pi \| \theta)$, note that if $\mu$ and $\theta$ are equivalent, then properties of the logarithm yield
\[
D_{\rm KL}(\pi\| \theta) = \int \log\left(\frac{d\mu}{d\theta}\right) \,d\pi + D_{\rm KL}(\pi\| \mu),
\]
and thus \eqref{eq:IP} fits the presented setting \eqref{primalProblem} using $c = \log\left(\frac{d\mu}{d\theta}\right)$.




Instead of solving \eqref{primalProblem} directly, we consider a dual formulation. We set
\begin{equation*}
\otimes_{i = 1}^N h_i \coloneqq \parentheses{h_1, \ldots, h_N} \in \prod\limits_{i = 1}^N H_i \eqqcolon H \ ,
\end{equation*}
to denote the respective elements in $ H $, the Cartesian product of the spaces $ H_i $. Similarly, given $ h \coloneqq \otimes_{j = 1}^N h_j $ and $i \in \{1, \dots, N\}$,
we will write
\begin{equation*}
    h_{-i} \coloneqq \otimes_{j \neq i} h_j \in \prod\limits_{j \neq i} H_j \eqqcolon H_{-i} \ ,
\end{equation*} 
to describe the tuple obtained from $ h $ by removing its $ i $-th component. In the same spirit, $ h_{-i} \otimes \tilde{h}_i \in H $
stands for the tuple obtained by replacing the $i$-th component of $ h $ by $ \tilde{h}_i \in H_i $. We also write 
\begin{equation*}
    \oplus_{i = 1}^N h_i \coloneqq \sum\limits_{i = 1}^N h_i, \quad \bigoplus\limits_{i = 1}^N H_i \coloneqq \sum\limits_{i = 1}^N H_i \ ,
\end{equation*}
although the sum does not necessarily have to be direct. 
Finally, let $ \psi $ stand for the convex conjugate of $ \varphi $, which we assume to be continuously differentiable. Let
\begin{equation*}
    \begin{aligned}
        F : H &\to \extendedReals \\
        \otimes_{i = 1}^N h_i &\mapsto \integral{\psi\parentheses{\oplus_{i = 1}^N h_i - c} - 
        \oplus_{i = 1}^N h_i}{\mu}
    \end{aligned}
\end{equation*}
and we denote by 
\begin{equation}
    \label{simplifiedDual}
    \tag{DP}
    \underset{h \in H}{\inf} F(h)
\end{equation}
the (up to sign change) dual problem to \eqref{primalProblem}. While this is structurally simply standard convex duality, of course there are subtleties that can occur in view of the infinite-dimensional nature of the spaces involved. In many recent works, strong duality is established in similar settings (see, e.g., \cite[Theorem A.1 and 2.2]{eckstein2021computation} or \cite{zaev2015monge}), with the main underlying assumption always being a kind of tightness of $\mathcal{Q}$. Since it is not the purpose of this paper to delve into the details of duality, we will simply assume throughout that strong duality holds (meaning the values of \eqref{primalProblem} and \eqref{simplifiedDual} correspond up to sign change), which we implicitly assume to include that the primal and dual versions of the IPFP correspond to each other as well (which essentially means we assume strong duality for each of the subspaces $H_i$ individually). That is, from now on we purely work on the dual side. It is worth mentioning that usually (cf.~\cite[Theorem 2.2]{eckstein2021computation}) a minimizer $\pi^*$ of \eqref{primalProblem} may be obtained from a minimizer $h^*$ of \eqref{simplifiedDual} via the formula
\[
\frac{d\pi^*}{d\mu} = \psi'\left(\oplus_{i=1}^N h_i^* - c\right),
\]
and thus applying the dual IPFP is not restrictive, even in terms of obtaining (approximate) optimizers for \eqref{primalProblem}.

We approach \eqref{simplifiedDual} via the dual version of the iterative proportional fitting procedure (IPFP), which is defined next.
Choose an arbitrary $ \CDIterate{h}{0} \in H $ and denote more generally its $ t $-th iterate by $ \CDIterate{h}{t} $, for $ t \in \integers_0 $.
The iterate $ \CDIterate{h}{t + 1} $ is obtained from its predecessor $ \CDIterate{h}{t} $ through intermediate iterates 
$ \CDIntermediateIterate{h}{t}{i} $, where $ i = 0, \ldots, N $, which in turn are inductively constructed as follows. 
Set $ \CDIntermediateIterate{h}{t}{0} \coloneqq \CDIterate{h}{t} $. To obtain $ \CDIntermediateIterate{h}{t}{i + 1} $ from 
$ \CDIntermediateIterate{h}{t}{i} $, select
\begin{equation}
    \label{coordinateSelection}
    \CDIntermediateIterate{h}{t}{i}_i \in \argmin\limits_{\tilde{h}_i \in H_i} F(\CDIntermediateIterate{h}{t}{i}_{-i} \otimes \tilde{h}_{-i}) \ ,
\end{equation}
and set $ \CDIntermediateIterate{h}{t}{i + 1} \coloneqq \CDIntermediateIterate{h}{t}{i}_{-i} \otimes \CDIntermediateIterate{h}{t}{i}_i $.
Finally, define $ \CDIterate{h}{t + 1} \coloneqq \CDIntermediateIterate{h}{t}{N} $. Hereby, we emphasize that by the above duality restricted to each of the subspaces $H_i$, the step \eqref{coordinateSelection} corresponds to the projection onto $\mathcal{Q}_i$ in the primal version of the IPFP. Similar to the above, given strong duality we have $\frac{d\pi^{(t)}}{d\mu} = \psi'\left(\oplus_{i=1}^N h_i^{(t)} - c\right)$. 
Throughout, we will assume that a sequence of iterates $ \parentheses{\CDIterate{h}{t}}_{t = 0}^\infty $ generated by the IPFP with
respect to a given starting point $ \CDIterate{h}{0} \in H $ exists. To shortly recap, the other standing assumptions are non-degeneracy of the spaces $H_1, \dots, H_N$, continuous differentiability of $\psi$, and strong duality.

\section{Exponential convergence of the IPFP}
\label{convergenceSection}

We first define the convexity and smoothness properties pertaining to $ \psi $, which will be needed throughout the paper.
\begin{definition}
    A differentiable function $ \psi : \mathbb{R} \rightarrow \mathbb{R} $ is called
    $ \parentheses{\sigma_R}_{R > 0} $-\textbf{strongly convex on bounded sets}, if for any $ R > 0 $, we have $ \sigma_R > 0 $ and
    \begin{equation*}
        \psi(\tilde{s}) \geq \psi^{\prime}(s)\parentheses{\tilde{s} - s} + \frac{\sigma_R}{2} \parentheses{\tilde{s} - s}^2 \ ,
    \end{equation*}
    for all $ \tilde{s}, s \in \openInterval{-R}{R} $.
\end{definition}
\begin{definition}
    A function $ \psi : \mathbb{R} \rightarrow \mathbb{R} $ is said to be 
    $ \parentheses{L_R}_{R > 0} $-\textbf{Lipschitz-smooth on bounded sets} if $ \psi $ is differentiable
    with derivative $ \psi^{\prime} $, and $ \psi^{\prime} $ is $ L_R $-Lipschitz-continuous on the bounded interval $ \openInterval{-R}{R} $,
    for all $ R > 0 $.
\end{definition}

We now define the key geometric property pertaining to the spaces $ H_i $, which underpins our convergence analysis.
\begin{definition}
    Let $ \closure{H}_i $ denote the closure of $ H_i \subset \LpSpace{\mu}{\infty} $ in $ \LpSpace{\mu}{2} $, the space $ \LpSpace{\mu}{2} $
    being endowed with its usual norm $ \lTwoNorm{\cdot} $. The spaces $ H_i $ are said to satisfy the \textbf{closed sum property} if 
    $ \bigoplus\limits_{i = 1}^N \closure{H}_i $ is a closed subspace of $ \LpSpace{\mu}{2} $.
\end{definition}

In view of the definition above, we shall always endow $ \bigoplus\limits_{i = 1}^N \closure{H}_i $ with 
$ \lTwoNorm{\cdot} $.
Similarly, we consider $ \closure{H} $, the closure of $ H $ in $ \LpSpace{\mu}{2}^N $, the latter being endowed with the scalar product
\begin{equation*}
    \sumInnerProduct{\otimes_{i = 1}^N h_i}{\otimes_{i = 1}^N \tilde{h}_i} \coloneqq \sum\limits_{i = 1}^N h_i \tilde{h}_i \ ,
\end{equation*}
and corresponding norm $ \sumNorm{\cdot} $. Recall that the norm $ \sumNorm{\cdot} $ induces the usual product topology on $ \LpSpace{\mu}{2}^N $.
Thus, $ \closure{H} = \prod_{i = 1}^N \closure{H}_i $.

Following Carlier's convergence analysis of the IPFP for the multimarginal optimal transport problem \cite{carlier2022linear}, we try to establish 
that the local strong convexity and Lipschitz-smoothness of $ \psi $, in some sense, carry over to $ F $.
To make this analogy more precise, we define the kind of sets with respect to which local properties of $ F $ will be understood.
\begin{definition}
    Let $ \supNorm{\cdot} $ denote the usual norm associated with $ \LpSpace{\mu}{\infty} $. 
    We call a subset of $ H $ \textbf{uniformly bounded in summation around $ c $ with radius} $ R > 0 $, if for all of its elements 
    $ \otimes_{i = 1}^N h_i $, we have $ \supNorm{\oplus_{i = 1}^N h_i - c} < R $.
\end{definition}

\subsection{Main idea}
\label{convergenceMainIdeaSubsection}

Recall that $ H \subset \LpSpace{\mu}{\infty}^N $. Using that $ \psi^{\prime} $ is continuous, it is easily seen that $ F $ is Gâteaux
differentiable with respect to $ \sumNorm{\cdot} $. The Gâteaux differential at any point $ h = \otimes_{i = 1}^N h_i \in H $ is given by
\begin{equation*}
    \begin{aligned}
        \gateauxDerivative{F}(h, \cdot) : H &\to \reals \\
        r = \otimes_{i = 1}^N r_i &\mapsto \gateauxDerivative{F}(h, r) = 
        \lTwoInnerProduct{\psi^{\prime} \parentheses{\oplus_{i = 1}^N h_i - c} - 1}{ \oplus_{i = 1}^N r_i} \ .
    \end{aligned}
\end{equation*}
Now consider two elements $ \tilde{h} = \otimes_{i = 1}^N \tilde{h}_i $ and $ h = \otimes_{i = 1}^N h_i $, taken from a
subset of $ H $, which is uniformly bounded in summation around $ c $ with radius $ R > 0 $. By assumption, the function $ \psi $ is strongly
convex over $ \openInterval{-R}{R} $ with modulus of convexity $ \sigma_{R} > 0 $. Thus,
\begin{multline*}
    F(\tilde{h}) - F(h) \geq 
        \integral{\parentheses{\psi^{\prime}\parentheses{\oplus_{i = 1}^N h_i - c} - 1} \oplus_{i = 1}^N \parentheses{\tilde{h}_i - h_i}}{\mu} \\
        + \frac{\sigma_{R}}{2} \integral{\absoluteValue{\oplus_{i = 1}^N \parentheses{\tilde{h}_i - h_i}}^2}{\mu} \ ,
\end{multline*}
and rewriting the right hand side in terms of $ \gateauxDerivative{F}(h, \cdot) $ and $ \lTwoNorm{\cdot} $ gives
\begin{equation}
    \label{problematicInequality}
    F(\tilde{h}) - F(h) \geq \gateauxDerivative{F}(h, \tilde{h} - h) 
    + \frac{\sigma_{R}}{2} \lTwoNorm{\oplus_{i = 1}^N \parentheses{\tilde{h}_i - h_i}}^2 \ .
\end{equation}
Recall that the Gâteaux derivative has been taken with respect to $ \sumNorm{\cdot} $. Since our subsequent convergence analysis relies on 
establishing an inequality of Polyak-Lojasiewicz type, the norm
\begin{equation*}
    \lTwoNorm{\oplus_{i = 1}^N \parentheses{\tilde{h}_i - h_i}}^2
\end{equation*} 
in \eqref{problematicInequality} should be replaced with a strictly positive multiple of the $ \sumNorm{\cdot} $-distance between $ \tilde{h} $ 
and $ h $. Unfortunately, the sum of the spaces $ H_i $ is not necessarily direct. Therefore, even in cases in which the 
$ \sumNorm{\cdot} $-distance between $ \tilde{h} $ and $ h $ is non-zero, the norm
\begin{equation*}
    \lTwoNorm{\oplus_{i = 1}^N \parentheses{\tilde{h}_i - h_i}}^2 \ ,
\end{equation*}
and consequently $ F(\tilde{h}) - F(h) $ as well as $ \gateauxDerivative{F}(h, \tilde{h} - h) $ might all vanish. In other words, the norm 
$ \sumNorm{\cdot} $ is too strong to establish strong convexity, as it does not take into account that both $ F $ and its Gâteaux derivatives are functions of 
$ \oplus_{i = 1}^N $ rather than $ \otimes_{i = 1}^N $.

To remedy the above situation, our key idea is to embed $ H $ into $ \parentheses{\closure{H}, \sumNorm{\cdot}} $ and to quotient the latter
by the kernel of
\begin{equation*}
    \begin{aligned}
        S: \parentheses{\closure{H}, \sumNorm{\cdot}} &\to 
        \parentheses{\bigoplus\limits_{i = 1}^N \closure{H_i}, \lTwoNorm{\cdot}} \\
        \otimes_{i = 1}^N h_i &\mapsto \oplus_{i = 1}^N h_i  \ .
    \end{aligned}
\end{equation*}
We will denote the kernel of $ S $ by $ \ker S $. We can then consider the factored version of $ F $, on the linear subspace
\begin{equation*}
    \factoredDomain \coloneqq \set{\equivalenceClass{h} \in \closure{H} / \ker S : \oplus_{i = 1}^N h_i \in \LpSpace{\mu}{\infty}} \ ,
\end{equation*}
endowed with the quotient norm $\quotientNorm{\equivalenceClass{h}} \coloneqq \underset{r \in \ker S}{\inf} \sumNorm{h + r}$.
It is immediate that $ \quotientNorm{\equivalenceClass{h}} \leq \sumNorm{h} $, for any $ h \in \closure{H} $. 
Indeed, we will see in Theorem \ref{factored_strong_convexity:thm} that the quotient norm is sufficiently weak to establish a suitable notion of strong convexity for
$ \hat{F} $, provided the spaces $ H_i $ satisfy the closed sum property. Along the way, we will realize in Corollary 
\ref{equivalence_norms_corollary} that the quotient norm is actually equivalent to the norm
\begin{equation*}
    \closure{H} / \ker S \ni \equivalenceClass{h} \mapsto \lTwoNorm{S h} \ .
\end{equation*}
Although the two norms are equivalent, the norm $ \quotientNorm{\cdot} $ is better suited to capture the componentwise iterations of the IPFP. 
This will become more clear in the proof of Lemma \ref{improvement:lemma}.

As opposed to the strong convexity, the Lipschitz-smoothness of $ F $ on subsets which are bounded in summation easily follows from the
respective Lipschitz-smoothness of $ \psi $ on bounded intervals.

With a Polyak-Lojasiewicz estimate and Lipschitz-smoothness at our disposal, we obtain the exponential convergence of the IPFP via standard 
arguments, commonly encountered in the analysis of first-order methods.

\subsection{Results}\label{subsec:results}

Our first objective is to establish that the factored version of $ F $, given by
\begin{equation*}
        \hat{F} : \factoredDomain \to \reals,~
        \equivalenceClass{h} \mapsto F(h)
\end{equation*}
is \textquote{locally strongly convex}. To this end, the following lemma asserts that the considered quotient norm is sufficiently weak, 
provided the closed sum property holds.
\begin{lemma}
    \label{norm_comparison_lemma}
    If the spaces $ H_i $ satisfy the closed sum property, then there exists $ \delta > 0 $ such that
    \begin{equation}
        \label{norm_comparison:eq}
        \quotientNorm{\equivalenceClass{h}} \leq \frac{1}{\delta} \lTwoNorm{Sh} \coloneqq \frac{1}{\delta} \lTwoNorm{\oplus_{i = 1}^N h_i} \ ,
    \end{equation}
    for all $ h \coloneqq \otimes_{i = 1}^N h_i \in \closure{H} $.
    \begin{proof}
        On account of Jensen's inequality, we have
        \begin{equation*}
            \lTwoNorm{S h}^2 = N^2 \lTwoNorm{\sum\limits_{i = 1} \frac{1}{N} h_i}^2 \leq N \sumNorm{h}^2 \ ,
        \end{equation*}
        and taking the square root on both sides, gives
        \begin{equation}
            \label{sum_operator_norm_bound:eq}
            \lTwoNorm{S h} \leq \sqrt{N} \sumNorm{h} \ .
        \end{equation}
        It follows that $ S $ is a bounded, surjective linear operator. Now, consider the factored version of $ S $, given by
        \begin{equation*}
            \begin{aligned}
                \hat{S} : \parentheses{\closure{H} / \ker S, \quotientNorm{\cdot}} &\to \parentheses{\bigoplus\limits_{i = 1}^N \closure{H}_i, 
                \lTwoNorm{\cdot}} \\
                \equivalenceClass{h} &\mapsto \oplus_{i = 1}^N h_i \ .
            \end{aligned}
        \end{equation*}
        For any $ \tilde{h} \in \equivalenceClass{h} $, \eqref{sum_operator_norm_bound:eq} gives
        \begin{equation*}
            \lTwoNorm{\hat{S} \equivalenceClass{h}} = \lTwoNorm{S h} = \lTwoNorm{S \tilde{h}} \leq \sqrt{N} \sumNorm{\tilde{h}} \ ,
        \end{equation*}
        and we conclude that
        \begin{equation*}
            \lTwoNorm{\hat{S} \equivalenceClass{h}} \leq \underset{\tilde{h} \in \equivalenceClass{h}}{\inf} \sqrt{N} \sumNorm{\tilde{h}}
            = \sqrt{N} \quotientNorm{\equivalenceClass{h}} \ .
        \end{equation*}
        Consequently, the operator $ \hat{S} $ is a bounded isomorphism on the Banach space 
        $ \parentheses{\closure{H} / \ker S, \quotientNorm{\cdot}} $, with
        \begin{equation}
            \label{factored_sum_operator_norm:eq}
            \lTwoNorm{\hat{S}\equivalenceClass{h}} \leq \sqrt{N} \quotientNorm{\equivalenceClass{h}} \ .
        \end{equation}
        By assumption, the range of $ \hat{S} $ is closed in $ \LpSpace{\mu}{2} $, and therefore a Banach space, when endowed
        with the relative norm induced by $ \lTwoNorm{\cdot} $. By virtue of the 
        Bounded Inverse Theorem \cite[Corollary 2.12 c)]{rudin1991functional} there exists $ \delta > 0 $ such that
        \begin{equation*}
            \quotientNorm{\equivalenceClass{h}} \leq \frac{1}{\delta} \lTwoNorm{\hat{S}\equivalenceClass{h}} \coloneqq
             \frac{1}{\delta} \lTwoNorm{\oplus_{i = 1}^N h_i} \ ,
        \end{equation*}
        and the proof is complete.
    \end{proof}
\end{lemma}

\begin{remark}
    Let $ \fixedSizeNorm{\hat{S}} $ denote the operator norm of $ \hat{S} $. Note that $\|S\| = \|\hat{S}\|$, which follows from the fact that $\|\hat{S}[h]\|_{L^2(\mu)} = \|Sh\|_{L^2(\mu)}$ for all $h \in [h]$. Further, by definition of the operator norm, \eqref{factored_sum_operator_norm:eq}
    is equivalent to $ \fixedSizeNorm{\hat{S}} \leq \sqrt{N} $. This upper bound might be very loose, depending on the considered spaces
    $ H_i $.
    Compare for instance with Example \ref{ex:comparisoncarlier}, where $\|\hat{S}\| = 1$ for all $N$.
\end{remark}

\begin{remark}
    \label{delta_operator_norm:remark}
    Let $ \hat{S}^{-1} $ denote the inverse of $ \hat{S} $, and $ \fixedSizeNorm{\hat{S}^{-1}} $ the respective operator norm. By construction, 
    the largest possible choice for $ \delta $ is given by $ 1 / \fixedSizeNorm{\hat{S}^{-1}} $.
\end{remark}

\begin{remark}\label{rem:iffinvertible}
    The closed sum property is both sufficient and necessary for \eqref{norm_comparison:eq} to hold. In fact, if \eqref{norm_comparison:eq} holds
    for some $ \delta > 0 $, then $ \hat{S} $ is continuously invertible. Consequently, the map
    \begin{equation*}
        \hat{S} : \parentheses{\closure{H} / \ker S, \quotientNorm{\cdot}} \to \parentheses{\bigoplus\limits_{i = 1}^N \closure{H}_i, 
        \lTwoNorm{\cdot}}
    \end{equation*}
    defines a homeomorphism between its domain and its codomain. Its domain being complete, its codomain must also be complete. Thus,
    $ \bigoplus\limits_{i = 1}^N \closure{H}_i $ must be a closed subspace of $ \LpSpace{\mu}{2} $.
\end{remark}

\begin{corollary}
    \label{equivalence_norms_corollary}
    If the spaces $ H_i $ satisfy the closed sum property, then \eqref{norm_comparison:eq} and \eqref{factored_sum_operator_norm:eq} show that $ \quotientNorm{\cdot} $ is equivalent to the norm $\closure{H} / \ker S \ni \equivalenceClass{h} \mapsto \lTwoNorm{S h}$.
\end{corollary}
Despite the above equivalence, we will continue to work with $ \quotientNorm{\cdot} $, the reasons of which will become more clear in
Lemma \ref{improvement:lemma}.
The Gâteaux differentiability of $ \hat{F} $ with respect to $ \quotientNorm{\cdot} $ is established in the next lemma.
\begin{lemma}
    For any $ \equivalenceClass{h} \in \factoredDomain $, the linear map
    \begin{equation*}
        \begin{aligned}
            \gateauxDerivative{\hat{F}}(\equivalenceClass{h}, \cdot) : \factoredDomain &\to \reals \\
            \equivalenceClass{r} &\mapsto d F(h, r) \ ,
        \end{aligned}
    \end{equation*}
    is bounded with respect to $ \quotientNorm{\cdot} $, and we have
    \begin{equation}
        \label{operator_norm_gateaux_derivative:eq}
        \quotientNorm{\gateauxDerivative{\hat{F}}(\equivalenceClass{h}, \cdot)} 
        \leq \sumNorm{\gateauxDerivative{F}(h, \cdot)} \ .
    \end{equation}
    for the respective operator norms.
    Moreover, the map $ \gateauxDerivative{\hat{F}}(\equivalenceClass{h}, \cdot) $ is the Gâteaux differential of $ \hat{F} $ at 
    $ \equivalenceClass{h} $.
    \begin{proof}
        To derive the first part, choose any $ \equivalenceClass{r} = \equivalenceClass{\otimes_{i = 1}^N r_i} \in H $. For all
        $ \tilde{r} \in \equivalenceClass{r} $ we have
        \begin{equation*}
            \absoluteValue{\gateauxDerivative{\hat{F}}(\equivalenceClass{h}, \equivalenceClass{r})}
            = \absoluteValue{\gateauxDerivative{F}(h, \tilde{r})} 
            \leq \sumNorm{\gateauxDerivative{F}(h, \cdot)} \sumNorm{\tilde{r}} \ .
        \end{equation*}
        Consequently,
        \begin{equation*}
            \absoluteValue{\gateauxDerivative{\hat{F}}(\equivalenceClass{h}, \equivalenceClass{r})} 
            \leq \sumNorm{\gateauxDerivative{F}(h, \cdot)} \underset{\tilde{r} \in \equivalenceClass{r}}{\inf} \sumNorm{\tilde{r}}
            = \sumNorm{\gateauxDerivative{F}(h, \cdot)} \quotientNorm{\equivalenceClass{r}} \ ,
        \end{equation*}
        whence
        \begin{equation*}
            \quotientNorm{\gateauxDerivative{\hat{F}}(\equivalenceClass{h}, \cdot)} \leq \sumNorm{\gateauxDerivative{F}(h, \cdot)} \ .
        \end{equation*}

        The second part is immediate in view of
        \begin{equation*}
            \frac{\hat{F}(\equivalenceClass{h + r}) - \hat{F}(\equivalenceClass{h}) 
            - \gateauxDerivative{\hat{F}}(\equivalenceClass{h}, \equivalenceClass{r})}{t}
            = \frac{F(h + r) - F(h) 
            - \gateauxDerivative{F}(h, r)}{t} \ ,
        \end{equation*}
        for any $ t \neq 0 $.
    \end{proof}
\end{lemma}

\begin{theorem}
    \label{factored_strong_convexity:thm}
    Suppose that the spaces $ H_i $ satisfy the closed sum property, and let $ \delta $ as in Lemma \ref{norm_comparison_lemma}. If, moreover,
    $ \psi $ is $ \parentheses{\sigma_R}_{R > 0} $-strongly convex, then for any subset $ \mathcal{M} \subset H $, which is bounded in summation
    around $ c $ with radius $ R > 0 $, it holds that
    \begin{equation}
        \label{factored_strong_convexity:eq}
        \hat{F}(\equivalenceClass{\tilde{h}}) - \hat{F}(\equivalenceClass{h}) 
        \geq \gateauxDerivative{\hat{F}}(\equivalenceClass{h}, \equivalenceClass{\tilde{h} - h}) 
        + \frac{\delta^2 \sigma_{R}}{2} \quotientNorm{\equivalenceClass{\tilde{h} - h}}^2 \ ,
    \end{equation}
    for any $ \tilde{h} $, $ h \in \mathcal{M} $.
    \begin{proof}
        Substituting the factored counterparts of $ F $ and its Gâteaux differential at $ h $ into \eqref{problematicInequality}, we get
        \begin{equation*}
            \hat{F}(\equivalenceClass{\tilde{h}}) - \hat{F}(h) 
            \geq \gateauxDerivative{\hat{F}}(\equivalenceClass{h}, \equivalenceClass{\equivalenceClass{\tilde{h} - h}})
            + \frac{\sigma_{R}}{2} \lTwoNorm{\oplus_{i = 1}^N \parentheses{\tilde{h}_i - h_i}}^2 \ .
        \end{equation*}
        This gives \eqref{factored_strong_convexity:eq} when combined with \eqref{norm_comparison:eq} from Lemma \ref{norm_comparison_lemma}.
    \end{proof}
\end{theorem}
\begin{corollary}[Polyak-Lojasiewicz inequality]
    \label{pl:corollary}
    Under the assumptions of Theorem \ref{factored_strong_convexity:thm} we have
    \begin{equation}
        \label{pl:eq}
        \tag{PL}
        F(\tilde{h}) - F(h) \geq- \frac{1}{2 \delta^2 \sigma_{R}} \quotientNorm{\gateauxDerivative{\hat{F}}(\equivalenceClass{h}, \cdot)}^2 \ .
    \end{equation}
    \begin{proof}
        Replacing the function values of $ \hat{F} $ with their corresponding function values with respect to $ F $, we conclude from 
        \eqref{factored_strong_convexity:eq} that
        \begin{equation*}
            F(\tilde{h}) - F(h) \geq \gateauxDerivative{\hat{F}}(\equivalenceClass{h}, \equivalenceClass{\tilde{h} - h}) 
            + \frac{\delta^2 \sigma_{R}}{2} \quotientNorm{\equivalenceClass{\tilde{h} - h}}^2 \ .
        \end{equation*}
        By the definition of the operator norm, we have
        \begin{equation*}
            \gateauxDerivative{\hat{F}}(\equivalenceClass{h}, \equivalenceClass{\tilde{h} - h})
            \geq - \absoluteValue{\gateauxDerivative{\hat{F}}(\equivalenceClass{h}, \equivalenceClass{\tilde{h} - h})}
            \geq - \quotientNorm{\gateauxDerivative{\hat{F}}(\equivalenceClass{h}, \cdot)} \quotientNorm{\equivalenceClass{\tilde{h} - h}} \ .
        \end{equation*}
        Using Young's inequality, we obtain
        \begin{equation*}
            \begin{split}
                - \quotientNorm{\gateauxDerivative{\hat{F}}(\equivalenceClass{h}, \cdot)} \quotientNorm{\equivalenceClass{\tilde{h} - h}}
            &= - \delta \sqrt{\sigma_{R}} \quotientNorm{\equivalenceClass{\tilde{h} - h}}
            \frac{\quotientNorm{\gateauxDerivative{\hat{F}}(\equivalenceClass{h}, \cdot)}}{\delta \sqrt{\sigma_{R}}} \\
            &\geq - \frac{1}{2 \delta^2 \sigma_{R}} \quotientNorm{\gateauxDerivative{\hat{F}}(\equivalenceClass{h}, \cdot)}^2
            - \frac{\delta^2 \sigma_{R}}{2} \quotientNorm{\equivalenceClass{\tilde{h} - h}}^2 \ .
            \end{split}
        \end{equation*}
        Combining all of the above estimates gives
        \begin{multline*}
            F(\tilde{h}) - F(h) \geq - \frac{1}{2 \delta^2 \sigma_{R}} \quotientNorm{\gateauxDerivative{\hat{F}}(\equivalenceClass{h}, \cdot)}^2
                - \frac{\delta^2 \sigma_{R}}{2} \quotientNorm{\equivalenceClass{\tilde{h} - h}}^2 \\
                + \frac{\delta^2 \sigma_{R}}{2} \quotientNorm{\equivalenceClass{\tilde{h} - h}}^2 \quad .
        \end{multline*}
        But the second and the last summand cancel, and \eqref{pl:eq} is proved.
    \end{proof}
\end{corollary}
Loosely speaking, \eqref{pl:eq} says that any $ h \in H $ for which the operator norm of $ \gateauxDerivative{\hat{F}}(\equivalenceClass{h}, \cdot) $ is small, the respective value $ F(h) $ cannot be much larger than any other
value $ F(\tilde{h}) $. 
Having established \eqref{pl:eq}, our next objective is to show that the IPFP makes large improvements
whenever
\begin{equation*}
    \quotientNorm{\gateauxDerivative{\hat{F}}(\equivalenceClass{\CDIterate{h}{t}})} \gg 1 \ 
\end{equation*}
at the given iterate $ \CDIterate{h}{t} $. By abuse of notation, we continue to write $ r_i $ and $ \equivalenceClass{r_i} $ for 
$  0_{-i} \otimes r_i $ and $ \equivalenceClass{ 0_{-i} \otimes r_i} $, respectively, given $ r_i \in H_i $. Now, consider the partial derivatives
\begin{equation*}
    \begin{aligned}
        \partialGateauxDerivative{F}{i}(h, \cdot) : H_i &\to \reals \\
        r_i &\mapsto \gateauxDerivative{F}(h, r_i) = \integral{\parentheses{\psi^{\prime}(\oplus_{i = 1}^N h_i - c) - 1} r_i}{\mu} \ ,
    \end{aligned}
\end{equation*}
at any $ h \in H $ and note that
\begin{equation*}
    \gateauxDerivative{F}(h, r) = \sum\limits_{i = 1}^N \partialGateauxDerivative{F}{i}(h, r_i) \ ,
\end{equation*}
for any $ r = \otimes_{i = 1}^N r_i $. We show that the maps $ h \mapsto \partialGateauxDerivative{F}{i}(h, \cdot) $ are Lipschitz-continuous
with respect to $ \sumNorm{\cdot} $.
\begin{lemma}
    \label{lipschitz:lemma}
    If $ \psi $ is $ \parentheses{L_R}_{R > 0} $-smooth and if $ \mathcal{M} \subset H $ is a subset which is bounded in summation around $ c $ 
    with radius $ R > 0 $, then
    \begin{equation*}
        \lTwoNorm{\partialGateauxDerivative{F}{i}(\tilde{h}, \cdot) - \partialGateauxDerivative{F}{i}(h, \cdot)}
        \leq L_{R} \fixedSizeNorm{S} \sumNorm{\tilde{h} - h} \ ,
    \end{equation*}
    for any $ i = 1, \ldots, N $ and $ \tilde{h}, h \in \mathcal{M} $.
    \begin{proof}
        For any $ r_i \in H_i $, we have
        \begin{equation*}
            \absoluteValue{\partialGateauxDerivative{F}{i}(\tilde{h}, r_i) - \partialGateauxDerivative{F}{i}(h, r_i)}
            \leq \integral{\absoluteValue{\psi^{\prime}(\oplus_{i = 1}^N \tilde{h}_i - c) - \psi^{\prime}(\oplus_{i = 1}^N h_i - c)}
            \absoluteValue{r_i}}{\mu} \ .
        \end{equation*}
        By assumption, $ \psi^{\prime} $ is Lipschitz-continuous over $ \openInterval{-R}{R} $ with Lipschitz-constant $ L_{R} $.
        It follows that,
        \begin{equation*}
            \absoluteValue{\partialGateauxDerivative{F}{i}(\tilde{h}, r_i) - \partialGateauxDerivative{F}{i}(h, r_i)}
            \leq L_{R} \integral{\absoluteValue{\oplus_{i = 1}^N (\tilde{h} - h)} \absoluteValue{r_i}}{\mu} \ .
        \end{equation*}
        Applying H\"{o}lder's inequality gives
        \begin{equation*}
            \begin{split}
                \absoluteValue{\partialGateauxDerivative{F}{i}(\tilde{h}, r_i) - \partialGateauxDerivative{F}{i}(h, r_i)}
                &\leq L_{R} \lTwoNorm{S(\tilde{h} - h)} \lTwoNorm{r_i}  \\
                &\leq L_{R} \fixedSizeNorm{S} \sumNorm{\tilde{h} - h} \lTwoNorm{r_i} \ ,
            \end{split}
        \end{equation*}
        which is the desired conclusion.
    \end{proof}
\end{lemma}

We are now in a position to lower bound the improvement made between $ F(\CDIterate{h}{t}) $ and its successor.
\begin{lemma}
    \label{improvement:lemma}
    If $ \psi $ is $ \parentheses{\sigma_R}_{R > 0} $-strongly convex and $ \parentheses{L_R}_{R > 0} $-Lipschitz smooth, and if the sequence of 
    iterates $ \parentheses{\CDIterate{h}{t}}_{t = 0}^\infty \subset H $ is contained in a subset which is bounded in summation around $ c $ with 
    radius $ R > 0 $, then
    \begin{align}
        \label{derivative_distance_iterates:eq}
        \sumNorm{\gateauxDerivative{F}(\CDIterate{h}{t}, \cdot)} &\leq L_{R} \sqrt{N} \fixedSizeNorm{S} 
        \sumNorm{\CDIterate{h}{t + 1} - \CDIterate{h}{t}} \\
        \label{distance_iterates_improvement_gap:eq}
        \sumNorm{\CDIterate{h}{t + 1} - \CDIterate{h}{t}}^2 
        &\leq \frac{2}{\sigma_{R}} \parentheses{F(\CDIterate{h}{t}) - F(\CDIterate{h}{t + 1})} \ .
    \end{align}
    Moreover,
    \begin{equation}
        \label{improvement_rate:eq}
        \frac{\sigma_{R}}{2 L_{R}^2 N \fixedSizeNorm{S}^2} \sumNorm{\gateauxDerivative{F}(\CDIterate{h}{t}, \cdot)}^2
        \leq F(\CDIterate{h}{t}) - F(\CDIterate{h}{t + 1}) \ .
    \end{equation}
    \begin{proof}
        We first derive \eqref{derivative_distance_iterates:eq}. Fix $ t \geq 0 $, and recall that
        \begin{equation*}
            \CDIntermediateIterate{h}{t}{i}_i \in 
            \argmin\limits_{\tilde{h}_i \in H_i} F(\CDIntermediateIterate{h}{t}{i - 1}_{-i} \otimes \tilde{h}_i) \ ,
        \end{equation*}
        by construction. It follows by first order conditions that
        \begin{equation}
            \label{fermat:eq}
            \partialGateauxDerivative{F}{i}(\CDIntermediateIterate{h}{t}{i}, \cdot) = 0 \ .
        \end{equation}
       	Consequently, we obtain for any $ r = \otimes_{i = 1}^N r_i \in H $ that
        \begin{equation*}
            \absoluteValue{\gateauxDerivative{F}(\CDIterate{h}{t}, r)} 
            = \absoluteValue{\sum\limits_{i = 1}^N \partialGateauxDerivative{F}{i}(\CDIterate{h}{t}, r_i)}
            \leq \sum\limits_{i = 1}^N \absoluteValue{\partialGateauxDerivative{F}{i}(\CDIntermediateIterate{h}{t}{i}, r_i)
            - \partialGateauxDerivative{F}{i}(\CDIterate{h}{t}, r_i)} \ .
        \end{equation*}
        Lemma \ref{lipschitz:lemma} now leads to
        \begin{equation*}
            \absoluteValue{\gateauxDerivative{F}(\CDIterate{h}{t}, r)} 
            \leq L_{R} \fixedSizeNorm{S} \sum\limits_{i = 1}^N \sumNorm{\CDIntermediateIterate{h}{t}{i} - \CDIterate{h}{t}}
            \lTwoNorm{r_i} \ .
        \end{equation*}
        But,
        \begin{equation*}
            \CDIntermediateIterate{h}{t}{i}_j
            = \begin{cases}
                \CDIterate{h}{t + 1}_j, &\text{for } j \leq i \\
                \CDIterate{h}{t}_j, &\text{otherwise} \ ,
            \end{cases}
        \end{equation*}
        and therefore
        \begin{equation*}
            \sumNorm{\CDIntermediateIterate{h}{t}{i} - \CDIterate{h}{t}} \leq \sumNorm{\CDIterate{h}{t + 1} - \CDIterate{h}{t}} \ .
        \end{equation*}
        We conclude that
        \begin{equation*}
            \absoluteValue{\gateauxDerivative{F}(\CDIterate{h}{t}, r)} 
            \leq L_{R} \fixedSizeNorm{S} \sumNorm{\CDIterate{h}{t + 1} - \CDIterate{h}{t}} \sum\limits_{i = 1}^N \lTwoNorm{r_i} \ ,
        \end{equation*}
        and H\"{o}lder's inequality for sums gives \eqref{derivative_distance_iterates:eq}.

        To obtain \eqref{distance_iterates_improvement_gap:eq}, begin by observing that $ \CDIntermediateIterate{h}{t}{i} $ and
        $ \CDIntermediateIterate{h}{t}{i - 1} $ only differ in the $ i $-th component by 
        $ \CDIntermediateIterate{h}{t}{i}_i - \CDIntermediateIterate{h}{t}{i - 1}_i $. Consequently, \eqref{problematicInequality} and
        \eqref{fermat:eq} immediately yield
        \begin{equation*}
            F(\CDIntermediateIterate{h}{t}{i - 1}) - F(\CDIntermediateIterate{h}{t}{i}) 
            \geq \frac{\sigma_{R}}{2} \lTwoNorm{\CDIntermediateIterate{h}{t}{i}_i - \CDIntermediateIterate{h}{t}{i-1}_i}^2 \ .
        \end{equation*}
        But,
        \begin{equation*}
            \CDIntermediateIterate{h}{t}{i}_i - \CDIntermediateIterate{h}{t}{i - 1}_i = \CDIterate{h}{t + 1}_i - \CDIterate{h}{t}_i \ .
        \end{equation*}
        Hence,
        \begin{equation*}
            F(\CDIntermediateIterate{h}{t}{i - 1}) - F(\CDIntermediateIterate{h}{t}{i}) 
            \geq \frac{\sigma_{R}}{2} \lTwoNorm{\CDIterate{h}{t + 1}_i - \CDIterate{h}{t}_i}^2 \ .
        \end{equation*}
        Summing the above estimate over $ i = 1, \ldots, N $ leads to \eqref{distance_iterates_improvement_gap:eq} after rearranging.

        As an immediate consequence of \eqref{derivative_distance_iterates:eq} of \eqref{distance_iterates_improvement_gap:eq} we obtain
        \eqref{improvement_rate:eq}, and the proof is complete.
    \end{proof}
\end{lemma}

\begin{remark}
    \label{norm_choice:remark}
    Employing the norm $ \sumNorm{\cdot} $ in the derivation of \eqref{derivative_distance_iterates:eq} introduced the attenuating factor $ N $ 
    in \eqref{derivative_distance_iterates:eq}, respectively $ 1 / N $ in \eqref{improvement_rate:eq}. But working with $ \sumNorm{\cdot} $  was 
    necessary for our proof of \eqref{distance_iterates_improvement_gap:eq}, which relied on the interplay $ \sumNorm{\cdot} $ and the 
    componentwise updates of the IPFP.
\end{remark}

\begin{corollary}
    If the sequence of iterates $ \parentheses{\CDIterate{h}{t}}_{t = 0}^\infty \subset H $ is contained in a subset which is bounded in
    summation around $ c $ with radius $ R > 0 $, we have
    \begin{equation}
        \label{factored_improvement_rate:eq}
        \frac{\sigma_{R}}{2 L_{R}^2 N \fixedSizeNorm{S}^2} 
        \quotientNorm{\gateauxDerivative{\hat{F}}(\equivalenceClass{\CDIterate{h}{t}}, \cdot)} 
        \leq F(\CDIterate{h}{t}) - F(\CDIterate{h}{t + 1}) \ .
    \end{equation}
    \begin{proof}
        Substituting \eqref{operator_norm_gateaux_derivative:eq} into \eqref{improvement_rate:eq} immediately gives 
        \eqref{factored_improvement_rate:eq}, as desired.
    \end{proof}
\end{corollary}

The exponential convergence of the IPFP is established by the next theorem.
\begin{theorem}
    \label{exponential_convergence:thm}
    Suppose the spaces $ H_i $ satisfy the closed sum property. If $ \psi $ is $ \parentheses{\sigma_R}_{R > 0} $-strongly convex and 
    $ \parentheses{L_R}_{R > 0} $-Lipschitz smooth, and if 
    the sequence of iterates $ \parentheses{\CDIterate{h}{t}}_{t = 0}^\infty \subset H $ is contained in a subset which is bounded in summation 
    around $ c $ with radius $ R > 0 $, also containing a global minimizer $ \overline{h} $, then
    \begin{equation}
        \label{exponential_rate:eq}
        F(\CDIterate{h}{t}) - F(\overline{h}) 
        \leq \parentheses{1 - \gamma}^t \parentheses{F(\CDIterate{h}{0}) - F(\overline{h})} \ ,
    \end{equation}
    for any $ t \geq 0 $, where
    \begin{equation}
        \label{rate_constant:eq}
        \gamma = \frac{1}{N} \parentheses{\frac{\sigma_{R}}{L_{R} \fixedSizeNorm{\hat{S}} \fixedSizeNorm{\hat{S}^{-1}}}}^2 \ .
    \end{equation}
    \begin{proof}
        By Lemma \ref{pl:corollary}, there exists $ \delta > 0 $ such that \eqref{pl:eq} holds. By Remark \ref{delta_operator_norm:remark}, 
        we may choose $ \delta = 1 / \fixedSizeNorm{\hat{S}^{-1}} $. Moreover, we have $ \fixedSizeNorm{S} = \fixedSizeNorm{\hat{S}} $. 
        Rearranging \eqref{pl:eq} and combining with \eqref{factored_improvement_rate:eq} gives
        \begin{equation}
            \label{relative_improvement:eq}
            \gamma
            \parentheses{F(\CDIterate{h}{t}) - F(\overline{h})} \leq F(\CDIterate{h}{t}) - F(\CDIterate{h}{t + 1}) \ .
        \end{equation}
        But,
        \begin{equation*}
            F(\CDIterate{h}{t}) - F(\CDIterate{h}{t + 1})
            = \parentheses{F(\CDIterate{h}{t}) - F(\overline{h})} + \parentheses{F(\overline{h}) - F(\CDIterate{h}{t + 1})} \ ,
        \end{equation*}
        and rearranging \eqref{relative_improvement:eq} yields
        \begin{equation*}
            F(\CDIterate{h}{t + 1}) - F(\overline{h}) 
        \leq \parentheses{1 - \gamma} \parentheses{F(\CDIterate{h}{t}) - F(\overline{h})} \ .
        \end{equation*}
        The claim now immediately follows by induction on $ t $.
    \end{proof}
\end{theorem}

\begin{example}\label{ex:comparisoncarlier}
	We can specify the bounds from Theorem \ref{exponential_convergence:thm} to the case of multi-marginal optimal transport, to showcase that the general result does not lose precision compared to the more specialized treatment in \cite{carlier2022linear}. In this case, $H_i$ are functions which only depend on the $i$-th variable of $\X = \X_1 \times \dots \times \X_N$, and in particular $H_1, \dots, H_{N-1}$ are normalized to have mean zero. Using the bounds from \cite[Lemma 3.1]{carlier2022linear} leads to $\|Sh-c\|_\infty \leq (4N-2) \|c\|_\infty$, that is, $L = \exp((4N-2)\|c\|_\infty)$ and $\sigma = \exp(-(4N-2))$, hence $(\sigma/L)^2 = \exp(-(16N-8)\|c\|_\infty)$ and thus the same constant as in \cite[Theorem 3.3]{carlier2022linear}. Hereby, as the spaces $H_i$ are orthogonal when using $\mu = \mu_1 \otimes \dots \otimes \mu_N$, we have \mbox{$\|\hat{S}\| = \|\hat{S}^{-1}\| = 1$.}
\end{example}

\section{Angles between subspaces and how they affect the rate of convergence}
\label{angles_section}

\subsection{General considerations}
\label{general_considerations_subsection}

In view of \eqref{rate_constant:eq}, we see that our rate of convergence deteriorates as 
$ \fixedSizeNorm{\hat{S}} \fixedSizeNorm{\hat{S}^{-1}} $, the condition number of $ \hat{S} $, becomes larger. It is therefore
desirable to derive respective upper bounds for $ \fixedSizeNorm{\hat{S}} $ and $ \fixedSizeNorm{\hat{S}^{-1}} $.

\subsubsection{Two subspaces}
The aim of this section is to derive upper bounds for $ \fixedSizeNorm{\hat{S}} $ and $ \fixedSizeNorm{\hat{S}^{-1}} $, when only
two subspaces are involved, i.e. when $ N = 2 $. We will see that the aforementioned operator norms relate to Friedrich's angle between 
$ \closure{H}_1 $ and $ \closure{H}_2 $.
\begin{definition}
    \label{friedrichsAngleDefinition}
    Denote by $ \closure{B}_1 $ the closed unit ball in $ \LpSpace{\mu}{2} $. Set
    \begin{equation*}
        \mathfrak{M} \coloneqq \closure{H}_1 \cap \closure{H}_2, \ \mathfrak{M}_1 \coloneqq \closure{H}_1 \cap \othorgonal{\mathfrak{M}}, \
        \mathfrak{M}_2 \coloneqq \closure{H}_2 \cap \othorgonal{\mathfrak{M}} \ ,
    \end{equation*}
    where $ \othorgonal{\mathfrak{M}} $ denotes the orthogonal complement of $ \mathfrak{M} $ in $ \LpSpace{\mu}{2} $. The angle 
    $ \alpha\parentheses{\closure{H}_1, \closure{H}_2} $ between $ \closure{H}_1 $ and $ \closure{H}_2 $ in 
    the sense of Friedrich is given by
    \begin{equation*}
        \alpha\parentheses{\closure{H}_1, \closure{H}_2} \coloneqq \arccos\parentheses{\friedrich{\closure{H}_1}{\closure{H}_2}} \in 
        \closedInterval{0}{\frac{\pi}{2}} \ ,
    \end{equation*}
    where
    \begin{equation*}
        \friedrich{\closure{H}_1}{\closure{H}_2} \coloneqq \sup \set{\absoluteValue{\innerProduct{h_1}{h_2}} : h_1 \in 
        \mathfrak{M}_1 \cap \closure{B}_1, \ h_2 \in \mathfrak{M}_2 \cap \closure{B}_1} \ .
    \end{equation*}
\end{definition}

\begin{remark}
    \label{closednessAngleEquivalenceRemark}
    In \cite{deutsch1995angle}, Deutsch showed that the spaces $ \closure{H}_1 $ and $ \closure{H}_2 $ satisfy the closed sum property if
    and only if $ \friedrich{\closure{H}_1}{\closure{H}_2} < 1 $.
\end{remark}

\begin{remark}
    \label{trivial_complemented_sum:remark}
    The space $ \mathfrak{M}_1 $ is trivial if and only if $ \closure{H}_1 \subset \closure{H}_2 $. Indeed, observe that
    $ \mathfrak{M}_1 $ is the space obtained by taking the orthogonal complement of $ \mathfrak{M} $ with respect to $ \closure{H}_1 $ instead
    of $ \LpSpace{\mu}{2} $. It follows that $ \closure{H}_1 = \mathfrak{M} \subset \closure{H}_2 $, whenever $ \mathfrak{M}_1 $ is trivial.
    Naturally, the same reasoning applies to $ \mathfrak{M}_2 $.
\end{remark}

\begin{theorem}
    \label{two_subspaces_bounds:thm}
    Suppose that $ H_1 $ and $ H_2 $ satisfy the closed sum property. Then,
    \begin{align}
        \label{operator_norm:eq}
        \fixedSizeNorm{\hat{S}} &= \begin{cases}
            \parentheses{1 + \friedrich{\closure{H}_1}{\closure{H}_2}}^{\frac{1}{2}} \ &\text{, if } \closure{H}_1 \cap \closure{H}_2 
            \text{ is trivial}, \\
            \sqrt{2} \ &\text{, otherwise},
        \end{cases} \\
        \label{norm_inverse:eq}
        \fixedSizeNorm{\hat{S}^{-1}} &= \parentheses{1 - \friedrich{\closure{H}_1}{\closure{H}_2}}^{-\frac{1}{2}} \ . &
    \end{align}
    \begin{proof}
        Using the orthogonal decompositions of $ \closure{H}_1 $ and $ \closure{H}_2 $ with respect to $ \mathfrak{M} $, we see that
        $ \closure{H}_1 + \closure{H}_2 $  may be decomposed as the direct sum of the subspaces $ \mathfrak{M} $, $ \mathfrak{M}_1 $ and
        $ \mathfrak{M}_2 $. Let 
        $ h_1 \in \closure{H}_1 $, $ h_2 \in \closure{H}_2 $ and denote by $ \tilde{h}_0 \in \mathfrak{M} $, $ \tilde{h}_1 \in \mathfrak{M}_1 $,
        $ \tilde{h}_2 \in \mathfrak{M}_2 $ the decomposition of $ h_1 + h_2 $ with respect to the subspaces $ \mathfrak{M} $, 
        $ \mathfrak{M}_1 $ and $ \mathfrak{M}_2 $. Because the decomposition is unique, and $ \ker S \subset \mathfrak{M}^2 $, while
        $ \mathfrak{M} $ is orthogonal to $ \mathfrak{M}_1 + \mathfrak{M}_2 $, we have
        \begin{equation}
            \begin{split}
                \label{explicit_quotient_norm:eq}
                \quotientNorm{\equivalenceClass{h_1 \otimes h_2}}^2 &= \sumNorm{ \parentheses{\tilde{h}_1 + \frac{1}{2} \tilde{h}_0} 
                \otimes \parentheses{\tilde{h}_2 + \frac{1}{2} \tilde{h}_0} }^2 \\
                &= \frac{1}{2} \lTwoNorm{\tilde{h}_0}^2 + \lTwoNorm{\tilde{h}_1}^2
                + \lTwoNorm{\tilde{h}_2}^2 \ .
            \end{split}
        \end{equation}
        On the other hand,
        \begin{equation*}
            \begin{split}
                \lTwoNorm{\hat{S} \equivalenceClass{h_1 + h_2}}^2 &= \lTwoNorm{\tilde{h}_0 + \tilde{h}_1 + \tilde{h}_2}^2 \\
                &= \lTwoNorm{\tilde{h}_0}^2 + \lTwoNorm{\tilde{h}_1 + \tilde{h}_2}^2 \\
                &= \lTwoNorm{\tilde{h}_0}^2 + \lTwoNorm{\tilde{h}_1}^2 + \lTwoNorm{\tilde{h}_2}^2 
                + 2 \lTwoInnerProduct{\tilde{h}_1}{\tilde{h}_2} \ .
            \end{split}
        \end{equation*}
        But
        \begin{equation*}
            \begin{split}
                2 \lTwoInnerProduct{\tilde{h}_1}{\tilde{h}_2} &\leq 2 \friedrich{\closure{H}_1}{\closure{H}_2} \lTwoNorm{\tilde{h}_1}
                \lTwoNorm{\tilde{h}_2} \\
                &\leq \friedrich{\closure{H}_1}{\closure{H}_2} \parentheses{\lTwoNorm{\tilde{h}_1}^2 + \lTwoNorm{\tilde{h}_2}^2} \ ,
            \end{split}
        \end{equation*}
        by definition of the respective Friedrich angle and Young's inequality. It follows that,
        \begin{multline}
            \label{estimate_factored_operator_norm:eq}
            \lTwoNorm{\hat{S} \equivalenceClass{h_1 + h_2}}^2 \leq
                \lTwoNorm{\tilde{h}_0}^2 \\
                + \parentheses{1 + \friedrich{\closure{H}_1}{\closure{H}_2}} \parentheses{\lTwoNorm{\tilde{h}_1}^2 + \lTwoNorm{\tilde{h}_2}^2} \ .
        \end{multline}
        Recall our standing assumption that neither $ \closure{H}_1 \subset \closure{H}_2 $ nor $ \closure{H}_2 \subset \closure{H}_2 $. By
        Remark \ref{trivial_complemented_sum:remark} this implies that neither $ \mathfrak{M}_1 $ nor $ \mathfrak{M}_2 $ is trivial.
        Now, if $ \mathfrak{M} \coloneqq \closure{H}_1 \cap \closure{H_2} $ is trivial, then $ \tilde{h}_0 $ is zero for any
        $ h_1 \in \closure{H}_1 $ and $ h_2 \in \closure{H}_2 $. Thus, we obtain
        \begin{equation*}
            \lTwoNorm{\hat{S} \equivalenceClass{h_1 + h_2}}^2 \leq \parentheses{1 + \friedrich{\closure{H}_1}{\closure{H}_2}}
            \quotientNorm{\equivalenceClass{h_1 \otimes h_2}}^2 \ ,
        \end{equation*}
        after substituting \eqref{explicit_quotient_norm:eq}. Using suitable sequences in $ \mathfrak{M}_1 $ and $ \mathfrak{M}_2 $ to approximate
        $ \friedrich{\closure{H}_1}{\closure{H}_2} $ via the respective inner products, we see that the above estimate is indeed sharp. 
        This proves the first case depicted in \eqref{operator_norm:eq}.
        If $ \mathfrak{M} $ is non-trivial, observe that $ \friedrich{\closure{H}_1}{\closure{H}_2} < 1 $, in view of the closed sum property
        and Remark \ref{closednessAngleEquivalenceRemark}. Hence, we 
        conclude from \eqref{estimate_factored_operator_norm:eq}
        that
        \begin{equation*}
            \begin{split}
                \lTwoNorm{\hat{S} \equivalenceClass{h_1 + h_2}}^2 &\leq 2 \parentheses{\frac{1}{2} \lTwoNorm{\tilde{h}_0}^2
                + \lTwoNorm{\tilde{h}_1}^2 + \lTwoNorm{\tilde{h}_2}^2} \\
                &= 2 \quotientNorm{\equivalenceClass{h_1 \otimes h_2}}^2 \ .
            \end{split}
        \end{equation*}
        Since the above estimate is sharp for any $ h_1 $, $ h_2 \in \mathfrak{M} $, the second case in \eqref{operator_norm:eq} is proved.

        Following the same reasoning, we see that
        \begin{multline*}
            \lTwoNorm{\hat{S} \equivalenceClass{h_1 \otimes h_2}}^2 \geq \lTwoNorm{\tilde{h}_0}^2 \\
            + \parentheses{1 - \friedrich{\closure{H}_1}{\closure{H}_2}} \parentheses{\lTwoNorm{\tilde{h}_1}^2
            + \lTwoNorm{\tilde{h}_2}^2} \ ,
        \end{multline*}
        the above estimate being sharp. Since $ 0 \leq \friedrich{\closure{H}_1}{\closure{H}_2} < 1 $, we have
        \begin{equation*}
            \begin{split}
                \parentheses{1 - \friedrich{\closure{H}_1}{\closure{H}_2}}^{-1} \lTwoNorm{\hat{S} \equivalenceClass{h_1 \otimes h_2}}^2
                &\geq \frac{1}{2} \lTwoNorm{\tilde{h}_0}^2 + \lTwoNorm{\tilde{h}_1}^2 + \lTwoNorm{\tilde{h}_2}^2 \\
                &= \quotientNorm{\equivalenceClass{h_1 \otimes h_2}}^2 \ ,
            \end{split}
        \end{equation*}
        no matter whether $ \mathfrak{M} $ is trivial or not. This gives \eqref{norm_inverse:eq}, and the proof is complete.
    \end{proof}
\end{theorem}

Recall that a linear map $P : \closure{H}_1 \times \closure{H}_2 \rightarrow \closure{H}_1 \times \closure{H}_2$ is called a (possibly oblique) projection onto $\ker(S)$ if $\Image(P) = \ker(S)$ and $P^2 = P$.
\begin{remark}\label{generalizedAngle:remark}
	Let $P$ be a projection onto $\ker(S)$. Then, it holds
    \begin{equation}
	\label{minimality_friedrich:eq}
	\friedrich{\closure{H}_1}{\closure{H}_2} \leq \sup \set{\frac{\absoluteValue{\lTwoInnerProduct{h_1}{h_2}}}{\lTwoNorm{h_1} \lTwoNorm{h_2}} : h_1 \otimes h_2 \in \ker(P)} .
\end{equation}
\begin{proof}
	Note that $ \specifiedNorm{\equivalenceClass{h}}{P} \coloneqq \sumNorm{(I - P) h} $ defines a norm over 
	$ \closure{H}_1 \times \closure{H}_2 / \ker(S) $. Moreover, choosing $ P $ to be the orthogonal projection onto $ \ker(S) $, recovers
	the quotient norm $ \quotientNorm{\cdot} $. It follows that $ \specifiedNorm{\cdot}{P} \geq \quotientNorm{\cdot} $, for any projection $ P $
	onto $ \ker(S) $, which implies $ \fixedSizeNorm{\hat{S}^{-1}}_P \geq \fixedSizeNorm{\hat{S}^{-1}} $, for the respective operator norms.
	Following the proof of Theorem \ref{two_subspaces_bounds:thm}, we see that $ \fixedSizeNorm{\hat{S}^{-1}}_P $ relates to
	\begin{equation*}
		\generalFriedrich{\closure{H}_1}{\closure{H}_2}{P} \coloneqq
		\sup \set{\frac{\absoluteValue{\lTwoInnerProduct{h_1}{h_2}}}{\lTwoNorm{h_1} \lTwoNorm{h_2}} : h_1 \otimes h_2 \in \ker(P)} \ ,
	\end{equation*}
	in the same way as $ \fixedSizeNorm{\hat{S}^{-1}} $ relates to $ \friedrich{\closure{H}_1}{\closure{H}_2} $, provided
	$ \generalFriedrich{\closure{H}_1}{\closure{H}_2}{P} < 1 $ (if it equals $1$, then the statement of the remark is trivially true). This gives \eqref{minimality_friedrich:eq}, as claimed.
\end{proof}
\end{remark}

\subsubsection{More than two subspaces}

We will now extend Theorem \ref{two_subspaces_bounds:thm} to the general case $ N \geq 2 $. To this end, we will
require the subspaces $ H_i $ to satisfy a strengthened version of the closed sum property.
\begin{definition}
    The subspaces $ H_i $, for $ i = 1, \ldots, N $, are said to satisfy the \textbf{strong closed sum property} if the spaces
    $ \bigoplus_{i = 1}^j \closure{H}_i $ are closed subspaces in $ \LpSpace{\mu}{2} $, for all $ j = 2, \ldots, N $.
\end{definition}

We extend the notation of the previous sections, writing $ \generalSumNorm{\cdot}{j} $, $ S_j $, $ \generalQuotientNorm{\cdot}{j} $ and
$ \hat{S}_j $, for $ j = 2, \ldots, N $, to denote the $ \ell^2 $-norm and the sum operator, alongside their factored versions which only 
involve the first $ j $ subspaces $ H_i $. Moreover, we define the sum operators
\begin{equation*}
    \begin{aligned}
        S_{j, j+1} : \parentheses{\bigoplus\limits_{i = 1}^j \closure{H}_i} \times \closure{H}_{j + 1} &\to 
        \bigoplus\limits_{i = 1}^{j + 1} \closure{H}_i \\
        \parentheses{\oplus_{i = 1}^j h_i} \otimes h_{j + 1} &\mapsto h_1 + \ldots + h_{j + 1} \ ,
    \end{aligned}
\end{equation*}
for $ j = 1, \ldots, N - 1 $, their domains being endowed with
\begin{equation*}
    \generalSumNorm{\parentheses{\oplus_{i = 1}^j h_i} \otimes h_{j + 1}}{j, j+1}^2 = \lTwoNorm{\oplus_{i = 1}^j h_i}^2
        + \lTwoNorm{h_{j + 1}}^2 \ ,
\end{equation*}
We denote the corresponding quotient norms with respect to $ \ker(S_{j, j + 1}) $ by $ \generalQuotientNorm{\cdot}{j, j + 1} $. As usual, we
denote by $ \hat{S}_{j, j + 1} $ the factored version of $ S_{j, j + 1} $ modulo $ \ker(S_{j, j + 1}) $, its domain being endowed with
$ \generalQuotientNorm{\cdot}{j, j + 1} $.
\begin{theorem}
    \label{general_subspaces:thm}
    If the spaces $ H_i $ satisfy the strong closed sum property, then
    \begin{equation}
        \label{norm_bound_N_subspaces:eq}
        \fixedSizeNorm{\hat{S}_N} \leq \prod\limits_{j = 1}^{N - 1} \fixedSizeNorm{\hat{S}_{j, j + 1}} \ ,
    \end{equation}
    and
    \begin{equation}
        \label{norm_bound_inverse_N_subspaces:eq}
        \fixedSizeNorm{\hat{S}_N^{-1}} \leq \prod\limits_{j = 1}^{N - 1} \fixedSizeNorm{\hat{S}_{j, j + 1}^{-1}} \ .
    \end{equation}
	Moreover, if the intersection of $\closure{H}_{\leq j} := \bigoplus\limits_{i = 1}^j \closure{H}_i $ with $ \closure{H}_{j + 1} $ is trivial for all
    $ j = 1, \ldots, N - 1 $, then
    \begin{equation}
        \label{condition_number_N_subspaces:eq}
        \parentheses{\fixedSizeNorm{\hat{S}_N} \fixedSizeNorm{\hat{S}_N^{-1}}}^2 \leq \prod\limits_{j = 1}^{N - 1}
        \left(\frac{1 + \friedrich{\closure{H}_{\leq j}}{\closure{H}_{j + 1}}}{1 - \friedrich{\closure{H}_{\leq j}}{\closure{H}_{j + 1}}}\right).
    \end{equation}
    \begin{proof}
        We only prove \eqref{norm_bound_inverse_N_subspaces:eq}, as \eqref{norm_bound_N_subspaces:eq} may be obtained in much the same way.
        The proof is by induction on $ N $. The base case $ N = 2 $ is immediate. Now assume the induction hypothesis holds for some $ N \geq 2 $.
        Consider $ h_i \in \closure{H}_i $, for $ i = 1, \ldots, N + 1 $. By definition of $ \hat{S}_{N, N + 1} $, we have
        \begin{equation}
            \label{NNplusOne:eq}
            \begin{split}
                \lTwoNorm{\oplus_{i = 1}^{N + 1} h_i}^2 &= \lTwoNorm{\parentheses{\oplus_{i = 1}^N h_i} \oplus h_{N + 1}}^2 \\
            &\geq \fixedSizeNorm{\hat{S}_{N, N + 1}^{-1}}^{-2} 
            \generalQuotientNorm{\equivalenceClass{\parentheses{\oplus_{i = 1}^N h_i} \otimes h_{N + 1}}}{N, N + 1}^2 \ .
            \end{split}
        \end{equation}
        By Hilbert's projection theorem, there exist $ r_i \in \closure{H}_i $, for $ i = 1, \ldots, N + 1 $ such that
        $ \parentheses{\oplus_{i = 1}^N r_i} \otimes r_{N + 1} \in \ker(S_{N, N + 1}) $ and
        \begin{equation*}
            \begin{split}
                \generalQuotientNorm{\equivalenceClass{\parentheses{\oplus_{i = 1}^N h_i} \otimes h_{N + 1}}}{N, N + 1}^2
            &= \generalSumNorm{\parentheses{\oplus_{i = 1}^N h_i - r_i} \otimes \parentheses{h_{N + 1} \oplus r_{N + 1}}}{N, N + 1}^2 \\
            &\coloneqq \lTwoNorm{\oplus_{i = 1}^N h_i - r_i}^2 + \lTwoNorm{h_{N + 1} \oplus r_{N + 1}}^2 \ .
            \end{split}
        \end{equation*}
        In turn, the induction hypothesis yields
        \begin{equation*}
            \lTwoNorm{\oplus_{i = 1}^N h_i - r_i}^2 \geq \prod\limits_{j = 1}^{N - 1} \fixedSizeNorm{\hat{S}_{j, j + 1}^{-1}}^{-2}
            \generalQuotientNorm{\equivalenceClass{\otimes_{i = 1}^N (h_i - r_i)}}{N}^2 \ .
        \end{equation*}
        Moreover, we immediately see from Theorem \ref{two_subspaces_bounds:thm} that $ \fixedSizeNorm{\hat{S}_{j, j + 1}^{-1}}^{-2} \leq 1 $,
        for all $ j = 1, \ldots, N - 1 $. It follows that,
        \begin{multline*}
            \generalQuotientNorm{\equivalenceClass{\parentheses{\oplus_{i = 1}^N h_i} \otimes h_{N + 1}}}{N, N + 1}^2 \\
            \geq \prod\limits_{j = 1}^{N - 1} \fixedSizeNorm{\hat{S}_{j, j + 1}^{-1}}^{-2} 
            \parentheses{\generalQuotientNorm{\equivalenceClass{\otimes_{i = 1}^N (h_i - r_i)}}{N}^2 + \lTwoNorm{h_{N + 1} \oplus r_{N + 1}}^2} \ .
        \end{multline*}
        Again, we conclude from Hilbert's projection theorem that there exist $ \tilde{r}_i \in \closure{H}_i $, for $ i = 1, \ldots, N $, such
        that $ \otimes_{i = 1}^N \tilde{r}_i \in \ker(S_N) $ and
        \begin{equation*}
            \generalQuotientNorm{\equivalenceClass{\otimes_{i = 1}^N (h_i - r_i)}}{N}^2 
            = \generalSumNorm{\otimes_{i = 1}^N (h_i - (r_i + \tilde{r}_i))}{N}^2 \ .
        \end{equation*}
        From this, we obtain
        \begin{multline*}
            \generalQuotientNorm{\equivalenceClass{\parentheses{\oplus_{i = 1}^N h_i} \otimes h_{N + 1}}}{N, N + 1}^2 \\
            \geq \prod\limits_{j = 1}^{N - 1} \fixedSizeNorm{\hat{S}_{j, j + 1}^{-1}}^{-2} 
            \parentheses{\generalSumNorm{\otimes_{i = 1}^N (h_i - (r_i \oplus \tilde{r}_i))}{N}^2 + \lTwoNorm{h_{N + 1} \oplus r_{N + 1}}^2} \ .
        \end{multline*}
        But $ \otimes_{i = 1}^N \tilde{r}_i \in \ker(S_N) $ and $ \parentheses{\oplus_{i = 1}^N r_i} \otimes r_{N + 1} \in \ker(S_{N, N + 1}) $.
        Consequently, we have $ \otimes_{i = 1}^N (r_i \oplus \tilde{r}_i) \otimes r_{N + 1} \in \ker(S_{N + 1}) $. Thus,
        \begin{equation*}
            \generalSumNorm{\otimes_{i = 1}^N (h_i - (r_i \oplus \tilde{r}_i))}{N}^2 + \lTwoNorm{h_{N + 1} \oplus r_{N + 1}}^2
            \geq \generalQuotientNorm{\equivalenceClass{\otimes_{i = 1}^{N + 1} h_i}}{N + 1}^2 \ .
        \end{equation*}
        Combining the two previous estimates yields
        \begin{equation}
            \label{norm_reduction:eq}
            \generalQuotientNorm{\equivalenceClass{\parentheses{\oplus_{i = 1}^N h_i} \otimes h_{N + 1}}}{N, N + 1}^2
            \geq \prod\limits_{j = 1}^{N - 1} \fixedSizeNorm{\hat{S}_{j, j + 1}^{-1}}^{-2}
            \generalQuotientNorm{\equivalenceClass{\otimes_{i = 1}^{N + 1} h_i}}{N + 1}^2 \ .
        \end{equation}
        Substituting \eqref{norm_reduction:eq} into \eqref{NNplusOne:eq} finally gives \eqref{norm_bound_inverse_N_subspaces:eq}.

        Combining \eqref{norm_bound_N_subspaces:eq} and \eqref{norm_bound_inverse_N_subspaces:eq}, and substituting
        $ \fixedSizeNorm{\hat{S}_{j, j + 1}} $ and $ \fixedSizeNorm{\hat{S}_{j, j + 1}^{-1}} $ with the corresponding values provided by
        Theorem \ref{two_subspaces_bounds:thm} immediately yields \eqref{condition_number_N_subspaces:eq}, which completes the proof.
    \end{proof}
\end{theorem}

\subsection{Martingale Optimal Transport}
\label{mot_subsection}

In this section, we derive an upper bound for the condition number $ \fixedSizeNorm{\hat{S}} \fixedSizeNorm{\hat{S}^{-1}} $
pertaining to the Martingale Optimal Transport problem (see, e.g., \cite{beiglbock2016problem, benamou2024entropic, chen2024convergence}). 
Let $ \mathcal{X}_1, \mathcal{X}_2 \subseteq \reals $ denote compact spaces, endowed with marginal distributions $ \mu_1 \in \mathcal{P}(\X_1) $, 
$ \mu_2 \in \mathcal{P}(\X_2) $, respectively. Set $ \mathcal{X} \coloneqq \mathcal{X}_1 \times \mathcal{X}_2 $ and assume a martingale coupling 
$ \mu $ between $ \mu_1 $ and $ \mu_2 $ exists. Within our framework, the marginal constraints may be formulated with respect to the function 
space
\begin{equation*}
    H_1 \coloneqq \set{f_1(x_1) + f_2(x_2) : f_1 \in \LpSpace{\mathcal{X}_1, \mu_1}{\infty}, \ f_2 \in \LpSpace{\mathcal{X}_2, \mu_2}{\infty}} 
    \subset \LpSpace{\mu}{\infty} \ ,
\end{equation*}
while the martingale constraint may be stated in terms of
\begin{equation*}
    H_2 \coloneqq \set{q(x_1)\Delta(x_1, x_2) : q \in \LpSpace{\mathcal{X}_1, \mu_1}{\infty}} \subset \LpSpace{\mu}{\infty} \  ,
\end{equation*}
where $ \Delta(x_1, x_2) \coloneqq \parentheses{x_2 - x_1} $.
Also recall that $ \mu_1 $ and $ \mu_2 $ must share the same mean, as $ \mu_1 $ is stochastically dominated by $ \mu_2 $, 
whenever a martingale coupling between $ \mu_1 $ and $ \mu_2 $ exists. Therefore, we may assume without loss of generality that $ \mu_1 $ and $ \mu_2 $ have zero mean.
We will adopt this assumption in the remainder of this section. We will also assume that $ \mu_2 $ is not a Dirac measure, for otherwise, the
problem becomes trivial. Finally, we assume that $ \mu $ is equivalent to $ \mu_1 \otimes \mu_2 $, the density
$ \density{\mu}{\parentheses{\mu_1 \otimes \mu_2}} $, being both essentially bounded and bounded away from zero.

We derive an upper bound for $ \fixedSizeNorm{\hat{S}} \fixedSizeNorm{\hat{S}^{-1}} $, using both, Theorem \ref{two_subspaces_bounds:thm}
and Remark \ref{generalizedAngle:remark}. By standard arguments, we obtain
\begin{align*}
    \closure{H}_1 &= \set{s(x_1, x_2) = f_1(x_1) + f_2(x_2) : f_1 \in \LpSpace{\reals, \mu_1}{2}, \ f_2 \in \LpSpace{\reals, \mu_2}{2}} \\
    \closure{H}_2 &= \set{q(x_1) \Delta(x_1, x_2) : q \in \LpSpace{\reals, \mu_1}{2}} \\
    \closure{H}_1 \cap \closure{H}_2 &= \set{a \Delta(x_1, x_2) : a \in \reals} \ .
\end{align*}
In particular, the closures of the spaces $ H_1 $ and $ H_2 $ in $ \LpSpace{\mu}{2} $ are the same for all $ \mu $ satisfying the above
assumptions. For ease of computation, we first ignore the fact that $\mu$ must be a martingale coupling and  assume $ \mu = \mu_1 \otimes \mu_2 $; the general case is obtained in Corollary \ref{cor:mot}.

\begin{theorem}
    \label{product_measure:thm}
    Set $V_2 \coloneqq \integral{x_2^2}{\mu_2}$ and denote by $ \absoluteValue{\mathcal{X}_1} $ the diameter of $ \mathcal{X}_1 $. If
    $ \mu = \mu_1 \otimes \mu_2 $, then
    \begin{equation}
        \label{angleMOT:eq}
        \friedrich{\closure{H}_1}{\closure{H}_2}
        \leq \frac{\absoluteValue{\mathcal{X}_1}}{\parentheses{\absoluteValue{\mathcal{X}_1}^2 + V_2}^{\frac{1}{2}}} \ .
    \end{equation}
    \begin{proof}
        Consider the oblique projection
        \begin{align*}
            P : \closure{H}_1 \times \closure{H}_2 &\to \closure{H}_1 \times \closure{H}_2 \\
            (f_1 + f_2) \otimes q \Delta &\mapsto - m_q \Delta \otimes m_q \Delta \ ,
        \end{align*}
        onto $ \ker(S) $, where $ m_q $ denotes the mean of $ q $ with respect to $ \mu_1 $. By Remark \ref{generalizedAngle:remark}, it suffices
        to prove that
        \begin{equation*}
            \sup \set{\frac{\absoluteValue{\lTwoInnerProduct{f_1 + f_2}{q \Delta}}}{\lTwoNorm{f_1 + f_2} \lTwoNorm{q \Delta}} : 
            (f_1 + f_2) \otimes q \Delta \in \ker(P)} \leq \frac{\absoluteValue{\mathcal{X}_1}}{\parentheses{\absoluteValue{\mathcal{X}_1}^2 + V_2}^{\frac{1}{2}}}.
        \end{equation*}
        Note that $ \ker P = \Image(I - P) = \closure{H}_1 \times H_2^0 $, where $ H_2^0 $ is the subspace of
        $ \closure{H}_2 $, comprised of the functions $ q \Delta $, for which $ m_q = 0 $. Now take $ (f_1 + f_2) \otimes q \Delta \in \ker P $, 
        and assume without loss of generality, that $ f_2 $ has zero mean, too. Since $ \mu = \mu_1 \otimes \mu_2 $, and $ f_2 $ as well as $ q $ 
        have zero mean, we obtain
        \begin{equation}
            \begin{split}
                \label{mot_inner_product:eq}
                \absoluteValue{\lTwoInnerProduct{f_1 + f_2}{q \Delta}} &= \absoluteValue{\integral{f_1(x_1) q(x_1) x_1 }{\mu_1(x_1)}} \\
                &\leq\specifiedNorm{f_1 q x_1}{\LpSpace{\mu_1}{1}} \\
                &\leq \specifiedNorm{f_1}{\LpSpace{\mu_1}{2}} \specifiedNorm{q x_1}{\LpSpace{\mu_1}{2}} \ ,
            \end{split}
        \end{equation}
        where the last estimate is due to H\"{o}lder's inequality. But
        \begin{equation*}
            \specifiedNorm{f_1}{\LpSpace{\mu_1}{2}} \leq \lTwoNorm{f_1 + f_2} \ ,
        \end{equation*}
        since $ f_2 $ has zero mean. Moreover,
        \begin{equation*}
            \specifiedNorm{q x_1}{\LpSpace{\mu_1}{2}}^2 = \integral{q^2 x_1^2}{\mu_1} \leq \absoluteValue{\mathcal{X}_1}^2 \integral{q^2}{\mu_1}
            = \frac{\absoluteValue{\mathcal{X}_1}^2}{V_2} \integral{q^2}{\mu_1} V_2 \ .
        \end{equation*}
        Consequently, we have
        \begin{equation*}
            \specifiedNorm{q x_1}{\LpSpace{\mu_1}{2}}^2 \leq (1 - \lambda) \integral{q^2 x_1^2}{\mu_1} 
            + \lambda \frac{\absoluteValue{\mathcal{X}_1}^2}{V_2} \integral{q^2}{\mu_1} V_2 \ ,
        \end{equation*}
        for all $ \lambda \in \openInterval{0}{1} $. Choosing $ \lambda $ such that
        \begin{equation*}
            (1 - \lambda) = \lambda \frac{\absoluteValue{\mathcal{X}_1}^2}{V_2} \ , \text{ that is } \lambda = \frac{V_2}{\absoluteValue{\mathcal{X}_1}^2 + V_2} \ ,
        \end{equation*}
        we obtain
        \begin{equation}
            \label{qDElta_norm:eq}
            \begin{split}
            \specifiedNorm{q x_1}{\LpSpace{\mu_1}{2}}^2 &\leq \frac{\absoluteValue{\mathcal{X}_1}^2}{\absoluteValue{\mathcal{X}_1}^2 + V_2}
            \parentheses{\integral{q^2 x_1^2}{\mu_1} + \integral{q^2}{\mu_1} V_2} \\
            &= \frac{\absoluteValue{\mathcal{X}_1}^2}{\absoluteValue{\mathcal{X}_1}^2 + V_2} \lTwoNorm{q \Delta}^2 \ ,
            \end{split}
        \end{equation}
        where the last equality follows again from our zero mean assumptions. Now, combining \eqref{mot_inner_product:eq} and 
        \eqref{qDElta_norm:eq} gives the desired conclusion.
    \end{proof}
\end{theorem}

In the following remark, we argue that the derived bound for the angle is even sharp in many cases.
\begin{remark}
Assume $\mu_1$ is symmetric. If there exists a sequence of symmetric densities 
$ \parentheses{E_n}_{n \in \mathbb{N}} $ such that $ E_n d\mu_1 $ converges weakly to 
$ \theta \coloneqq (\delta_{-\absoluteValue{\mathcal{X}_1}} + \delta_{\absoluteValue{\mathcal{X}_1}}) / 2 $, then estimate 
\eqref{angleMOT:eq} is sharp. In particular, this is the case whenever $\mu_1$ is equivalent to the uniform distribution on 
$[-|\mathcal{X}_1|, |\mathcal{X}_1|]$, or if $\theta$ is absolutely continuous with respect to $\mu_1$.
    \begin{proof}
        Denote by $ \sgn $ the sign function and set $ q_n(x_1) \coloneqq \sgn(x_1) \sqrt{E_n(x_1)} $, 
        $ f_{1, n}(x_1) \coloneqq \sqrt{E_n(x_1)} $ and $ f_{2, n} = 0 $. Note that
        \begin{equation*}
            f_{1, n} + f_{2, n} \in \closure{H}_1 \cap  \othorgonal{\parentheses{\closure{H}_1 \cap \closure{H}_2}}, \ 
            q_n \Delta \in \closure{H}_2 \cap \othorgonal{\parentheses{\closure{H}_1 \cap \closure{H}_2}} \ ,
        \end{equation*}
        by symmetry of $ E_n $ and $ \mu_1 $. Moreover,
        \begin{equation*}
            \lTwoInnerProduct{f_{1, n} + f_{2, n}}{q_n \Delta} = \integral{\absoluteValue{x_1} E_n(x_1)}{\mu_1(x_1)} 
                \to \integral{\absoluteValue{x_1}}{\theta}
                = \absoluteValue{\mathcal{X}_1} \ ,
        \end{equation*}
        as $ n \to \infty $. On the other hand, we have
        \begin{equation*}
            \lTwoNorm{q_n \Delta}^2 = \integral{q_n^2}{\mu_1} V_2 + \integral{q_n^2(x_1) x_1^2}{\mu_1} \to \absoluteValue{X_1}^2 + V_2
        \end{equation*}
        and $ \lTwoNorm{f_{1, n} + f_{2, n}} = 1 $, which yields
        \begin{equation*}
            \lTwoInnerProduct{\frac{f_{1, n} + f_{2, n}}{\lTwoNorm{f_{1, n} + f_{2, n}}}}{\frac{q_n \Delta}{\lTwoNorm{q_n \Delta}}} \to 
            \frac{\absoluteValue{\mathcal{X}_1}}{\parentheses{\mathcal{X}_1 + V_2}^{\frac{1}{2}}} \ ,
        \end{equation*}
        that is, the sharpness of \eqref{angleMOT:eq}.

        When $ \mu_1 $ is equivalent to the uniform distribution on $[-|\mathcal{X}_1|, |\mathcal{X}_1|] $ the densities $ E_n $ may be obtained
        by applying standard smoothing arguments to $ \theta $. When $ \theta $ is absolutely continuous with respect to $ \mu_1 $, we may
        obviously choose $ E_n = d\theta / d\mu_1 $, and the proof is complete.
    \end{proof}
\end{remark}

We are now in a position to bound the condition number for more general measures $ \mu $, as described in the beginning of this section.
The measure associated with the considered quotient and operator norms will be indicated with an additional subscript. E.g.\ we write 
$ \specifiedNorm{\cdot}{\sim, \mu} $ and $ \specifiedNorm{\cdot}{\sim, (\mu_1 \otimes \mu_2)} $ to distinguish between the quotient norms
associated with $ \mu $ and $ \mu_1 \otimes \mu_2 $, respectively.
\begin{corollary}\label{cor:mot}
    If $ \mu $ is equivalent to $ \mu_1 \otimes \mu_2 $, such that $ p \coloneqq \density{\mu}{\mu_1 \otimes \mu_2} $ is both, essentially bounded
    and essentially bounded away from zero, then
    \begin{align*}
        \fixedSizeNorm{\hat{S}}_{\sim, \mu} = \sqrt{2}, \quad
        \fixedSizeNorm{\hat{S}^{-1}}_{\sim, \mu} 
        \leq
        \sqrt{\frac{\supNorm{p^{-1}}}{\supNorm{p}}}
        \parentheses{1 - \frac{\absoluteValue{\mathcal{X}_1}}{\parentheses{\absoluteValue{\mathcal{X}_1}^2 + V_2}^{\frac{1}{2}}}}^{- \frac{1}{2}} \ .
    \end{align*}
    \begin{proof}
        The first equality follows immediately from Theorem \ref{two_subspaces_bounds:thm}. To obtain the second estimate, we first deduce from
        Theorem \ref{angleMOT:eq} and Theorem \ref{product_measure:thm} that
        \begin{equation*}
            \fixedSizeNorm{\hat{S}^{-1}}_{\sim, \mu_1 \otimes \mu_2} \leq
            \parentheses{1 - \frac{\absoluteValue{\mathcal{X}_1}}{\parentheses{\absoluteValue{\mathcal{X}_1}^2 + V_2}^{\frac{1}{2}}}}^{- \frac{1}{2}} \ .
        \end{equation*}
        Moreover, the assumptions on $ \mu $ yield
        \begin{equation*}
            \supNorm{p^{-1}}^{-\frac{1}{2}} \specifiedNorm{\cdot}{\LpSpace{\mu_1 \otimes \mu_2}{2}} 
            \leq \specifiedNorm{\cdot}{\LpSpace{\mu}{2}}
            \leq \supNorm{p}^{\frac{1}{2}} \specifiedNorm{\cdot}{\LpSpace{\mu_1 \otimes \mu_2}{2}} \ .
        \end{equation*}
        The above estimates remain true, when replacing $ \specifiedNorm{\cdot}{\LpSpace{\mu_1 \otimes \mu_2}{2}} $ and
        $ \specifiedNorm{\cdot}{\LpSpace{\mu}{2}} $ with the respective quotient norms $ \specifiedNorm{\cdot}{\sim, \mu_1 \otimes \mu_2} $
        and $ \specifiedNorm{\cdot}{\sim, \mu} $. Combining the upper bound on $ \fixedSizeNorm{\hat{S}^{-1}}_{\sim, \mu_1 \otimes \mu_2} $, with 
        the equivalence of the above norms yields the desired conclusion.
    \end{proof}
\end{corollary}

\section*{Acknowledgment}
The authors are grateful for support by the German Research Foundation through the project 553088969 as well as the Cluster of Excellence ``Machine Learning - New Perspectives for Science'' (EXC 2064/1 number 390727645).

\printbibliography
\end{document}